\documentclass[a4paper,11pt]{article}

\usepackage{amsmath}
\usepackage{amsthm}
\usepackage{amssymb}
\usepackage{cancel}
\usepackage{color}
\usepackage{calrsfs}
\usepackage{varioref}
\usepackage{bbm}
\definecolor{dblue}{rgb}{0.09,0.32,0.44} 
\usepackage[pdfborder={0 0 0}, colorlinks=true, pdfborderstyle={}, linkcolor=dblue, citecolor=dblue, urlcolor=blue]{hyperref}
\usepackage{epsfig}
\usepackage{psfrag}
\usepackage{subfigure}
\usepackage{graphicx}
\usepackage[mathscr,mathcal]{eucal}
\usepackage[shortlabels]{enumitem}
\usepackage{t1enc}
\usepackage[latin2]{inputenc}

\newtheorem {theorem}{Theorem}
\newtheorem {lemma}{Lemma}

\newtheorem {proposition}{Proposition}

\newtheorem* {theorem*}{Theorem}
\newtheorem* {lemma*}{Lemma}
\newtheorem* {corollary*}{Corollary}
\newtheorem* {proposition*}{Proposition}
\newtheorem* {definition*}{Definition}
\newtheorem* {conjecture*}{Conjecture}
\newtheorem* {question*}{Question}

\newtheorem* {theoremkv*} {Theorem KV}
\newtheorem* {corollarykv*} {Corollary KV}
\newtheorem* {theoremrsc1*} {Theorem RSC1}
\newtheorem* {theoremrsc2*} {Theorem RSC2}

\theoremstyle{remark}
\newtheorem* {remark*}{Remark}

\def \C {\mathbb C}

\def \N {\mathbb N}

\def \R {\mathbb R}

\def \Z {\mathbb Z}

\def\boP{\mathbf{P}}

\def\cB{\mathcal{B}}
\def\cC{\mathcal{C}}

\def\cE{\mathcal{E}}
\def\cF{\mathcal{F}}

\def\cH{\mathcal{H}}

\def\cL{\mathcal{L}}

\def\cN{\mathcal{N}}
\def\cO{\mathcal{O}}

\def\cU{\mathcal{U}}

\def\vareps{\varepsilon}

\newcommand{\probab}[1]{\ensuremath{\mathbf{P}\left(#1\right)}}
\newcommand{\expect}[1]{\ensuremath{\mathbf{E}\left(#1\right)}}
\newcommand{\var}[1]{\ensuremath{\mathbf{Var}\left(#1\right)}}
\newcommand{\cov}[2]{\ensuremath{\mathbf{Cov}\left(#1,#2\right)}}
\newcommand{\condprobab}[2]{\ensuremath{\mathbf{P}\left(#1\bigm|#2\right)}}
\newcommand{\condexpect}[2]{\ensuremath{\mathbf{E}\left(#1\bigm|#2\right)}}

\newcommand{\probabom}[1]{\ensuremath{\mathbf{P}_{\omega}\left(#1\right)}}
\newcommand{\expectom}[1]{\ensuremath{\mathbf{E}_{\omega}\left(#1\right)}}
\newcommand{\varom}[1]{\ensuremath{\mathbf{Var}_{\omega}\left(#1\right)}}

\newcommand{\condprobabom}[2]{\ensuremath{\mathbf{P}_{\omega}\left(#1\bigm|#2\right)}}
\newcommand{\condexpectom}[2]{\ensuremath{\mathbf{E}_{\omega}\left(#1\bigm|#2\right)}}

\newcommand{\ind}[1]{\ensuremath{\mathbbm{1}_{\{#1\}}}}

\def \toprob {\,\,\buildrel\boP\over\longrightarrow\,\,}

\def\clap#1{\hbox to 0pt{\hss#1\hss}}

\DeclareMathOperator*{\infsuplim}{\overline{\underline{\lim}}}

\def\sgn{\mathrm{sgn}}
\def \Ordo {\cO}
\def\ordo{o}

\DeclareMathOperator{\lap}{lap}
\DeclareMathOperator{\grad}{grad}
\let\div\undefined
\DeclareMathOperator{\div}{div}

\def\Ker{\mathrm{Ker}}
\def\Ran{\mathrm{Ran}}
\def\Dom{\mathrm{Dom}}

\def\sqd{\abs{\Delta}^{1/2}}
\def\nsqd{\abs{\Delta}^{-1/2}}

\renewcommand{\d}{\mathrm d}

\newcommand{\abs}[1]{\ensuremath\left|{#1}\right|}
\newcommand{\norm}[1]{\ensuremath\left\|{#1}\right\|}
\newcommand{\sprod}[2]{\ensuremath\left\langle{#1,#2}\right\rangle}

\def \wt {\widetilde}

\def\wh{\widehat}

\newcommand{\tostoptop}{\stackrel{\mathrm{st.op.top.}}{\longrightarrow}}

\title{Central Limit Theorem for Random Walks in Doubly Stochastic Random Environment: ${\cH_{-1}}$ Suffices}

\author{
{\sc Gady Kozma$^{1}$}
\qquad\qquad
{\sc B\'alint T\'oth$^{2,3}$}
\\[8pt]
{$^1$ Weizmann Institute, Rehovot, IL}
\\
{$^2$ School of Mathematics, University of Bristol, UK}
\\
{$^3$ R\'enyi Institute, Budapest, HU}
}

\begin{document}

\advance\baselineskip by 3pt
\maketitle

\begin{abstract}
We prove a central limit theorem under diffusive scaling for the displacement of a random walk on $\Z^d$ in stationary and ergodic doubly stochastic random environment, under the ${\cH_{-1}}$-condition imposed on the drift field. The condition is equivalent to  assuming that the stream tensor of the drift field be stationary and square integrable. This improves the best existing result \cite{komorowski_landim_olla_12}, where it is assumed that the stream tensor is in $\cL^{\max\{2+\delta, d\}}$, with $\delta>0$. Our proof relies on an extension of the \emph{relaxed sector condition} of \cite{horvath_toth_veto_12}, and is technically rather simpler than existing earlier proofs of similar results by Oelschl\"ager \cite{oelschlager_88} and Komorowski, Landim and Olla \cite{komorowski_landim_olla_12}.

\medskip\noindent
{\sc MSC2010: 60F05, 60G99, 60K37}

\medskip\noindent
{\sc Key words and phrases:} random walk in random environment, central limit theorem, Kipnis-Varadhan theory, sector condition.

\end{abstract}

\section{Introduction: setup and main result}
\label{s:Introduction: setup and main result}

Since its appearance in  the probability and physics literature in the mid-seventies the general topics of \emph{random walks/diffusions in random environment} became the most complex and robust area of research. For a general overview of the subject and its historical development we refer the reader to the surveys Kozlov \cite{kozlov_85}, Zeitouni \cite{zeitouni_01}, Biskup \cite{biskup_11} or Kumagai \cite{kumagai_14}, written at various stages of this rich story. The main problem considered in our paper is that of diffusive limit in the doubly stochastic (and hence, a priori stationary) case.

\subsection{The random walk and the \texorpdfstring{${\cH_{-1}}$}{H -1}-condition}
\label{ss:The random walk}

Let $(\Omega, \cF, \pi, \tau_z:z\in\Z^d)$ be a probability space with
an ergodic $\Z^d$-action. Denote by  $\cE_{+}:=\{e_1,\dots,e_d:
e_i\in\Z^{d}, \ \ e_i\cdot e_j=\delta_{i,j}\}$ the standard generating
basis in $\Z^d$ and let $\cE:=\{\pm e_j: e_j\in \cE_{+}\} =
\{k\in\Z^d: |k|=1\}$ be the set of possible steps of a
nearest-neighbour walk on $\Z^d$. Assume that bounded measurable
functions $p_k:\Omega\to[0,s^*]$, $k\in\cE$ are given ($s^*$ denotes the common bound), and assume the $p_k$ satisfy \emph{bistochasticity}, by which we mean the following property
\begin{align}
\label{bistoch}
\sum_{k\in\cE}p_k(\omega)
=
\sum_{k\in\cE}p_{-k}(\tau_k\omega).
\end{align}
Lift these functions to the lattice $\Z^d$ by defining
\begin{align}
\label{jump_probab_field}
P_k(x)=P_k(\omega,x):=p_k(\tau_x\omega).
\end{align}
Given these,  define the continuous time nearest neighbour random walk $X(t)$ as a continuous time Markov chain on $\Z^d$, with  $X(0)=0$ and conditional jump rates
\begin{align}
\label{the walk}
\condprobabom{X(t+dt)= x+k}{X(t)=x} = P_k(\omega, x) dt + \ordo(dt),
\end{align}
where the subscript $\omega$ denotes that the random walk $X(t)$ is a
continuous time Markov chain on $\Z^d$ \emph{conditionally}, with
$\omega\in\Omega$ sampled according to $\pi$. Note that \eqref{bistoch} is equivalent to 
\begin{align}
\notag
\sum_{k\in\cE}P_k(\omega,x)= \sum_{k\in\cE}P_{-k}(\omega,x+k),
\end{align}
which is exactly bistochasticity of the random walk defined in \eqref{the walk} above. Since the $p_k$ are bounded, so will be the  total jump rate of the walk
\begin{align}
\notag
p(\omega):=\sum_{k\in\cE}p_k(\omega)\le 2d s^*.
\end{align} 
Thus, there is no difference between the long time asymptotics of this walk and the discrete time (possibly lazy) walk $n\mapsto X_n\in\Z^d$ with jump probabilities 
\begin{align}
\notag
\condprobabom{X_{n+1}= y}{X_n=x} = 
\begin{cases}
(2d s^*)^{-1} P_k(\omega, x) 
&\text{if }y-x=k\in\cE, 
\\[2pt]
1-(2d s^*)^{-1} \sum_{l\in\cE}P_l(\omega, x)
&\text{if }y-x=0, 
\\[2pt]
0
&\text{if }y-x\not\in \cE \cup\{0\}.
\end{cases}
\end{align}
We speak about continuous time walk only for reasons of convenience, in order to easily quote facts and results form Kipnis-Varadhan theory of CLT for additive functionals of Markov processes, without tedious reformulations.

We formulate our problem and prove our main result in the context of nearest neighbour walks. This is only for convenience reason. The main result of this paper holds true for finite range bistochastic RWREs under the appropriate conditions. For more details on this see the remark after Theorem 1, further down in the paper.  

We will use the notation $\probabom{\cdot}$, $\expectom{\cdot}$  and $\varom{\cdot}$ for \emph{quenched} probability, expectation and variance. That is: probability, expectation, and variance with respect to the distribution of the random walk $X(t)$, \emph{conditionally, with given fixed environment $\omega$}. The notation $\probab{\cdot}:=\int_\Omega\probabom{\cdot} {\d}\pi(\omega)$, $\expect{\cdot}:=\int_\Omega\expectom{\cdot} {\d}\pi(\omega)$ and $\var{\cdot}:=\int_\Omega\varom{\cdot} {\d}\pi(\omega) + \int_\Omega\expectom{\cdot}^2 {\d}\pi(\omega) - \expect{\cdot}^2$ will be reserved for \emph{annealed} probability, expectation and variance. That is: probability,  expectation and variance with respect to the random walk trajectory $X(t)$ \emph{and} the environment $\omega$, sampled according to the distribution $\pi$. 

It is well known (and easy to check, see e.g.\ \cite{kozlov_85}) that due to double stochasticity \eqref{bistoch} the annealed set-up is stationary and ergodic in time: the process of \emph{the environment as seen from the position of the random walker} (to be formally defined soon) is stationary and ergodic in time under the probability measure $\pi$ and consequently the random walk $t\mapsto X(t)$ will have stationary and ergodic annealed increments.

Next we define, for $k\in\cE$, $v_k:\Omega\to[-s^*,s^*]$,
$s_k:\Omega\to[0,s^*]$,  and $\psi,\varphi:\Omega\to\R^d$,
\begin{align}
\label{v and phi}
&
v_k(\omega):=\frac{p_k(\omega)-p_{-k}(\tau_k\omega)}{2},
&&
\varphi(\omega):= \sum_{k\in\cE} k v_k(\omega),
\\
\label{s and psi}
&
s_k(\omega):=\frac{p_k(\omega)+p_{-k}(\tau_k\omega)}{2},
&&
\psi(\omega):= \sum_{k\in\cE} k s_k(\omega).
\end{align}
Their corresponding lifting to $\Z^d$ are
\begin{align}
\label{vlifting}
&
V_k(x)=V_k(\omega,x):=v_k(\tau_x\omega),
&&
\Phi(x)=\Phi(\omega,x):=\varphi(\tau_x\omega),
\\[8pt]
\notag
&
S_k(x)=S_k(\omega,x):=s_k(\tau_x\omega),
&&
\Psi(x)=\Psi(\omega,x):=\psi(\tau_x\omega).
\end{align}
Note that
\begin{align}
\notag
-s^*\le v_k(\omega)\le s^*,
\qquad
0\le s_k(\omega)\le s^*,
\qquad
\abs{\varphi(\omega)}\le 2\sqrt{d} s^*,
\qquad
\abs{\psi(\omega)}\le \sqrt{d} s^*,
\qquad
\text{a.s.}
\end{align}
The local \emph{quenched} drift of the random walk is
\begin{align}
\notag
\condexpectom{dX(t)}{X(s): 0\le s\le t} = \left(\Psi(\omega,X(t))+\Phi(\omega,X(t))\right) dt +\ordo(dt).
\end{align}
Note that from \eqref{bistoch} and the definitions \eqref{v and phi}, \eqref{s and psi} it follows that for $\pi$-almost all $\omega\in\Omega$
\begin{align}
\label{divfree}
v_k(\omega)&=-v_{-k}(\tau_k\omega),&
\sum_{k\in\cE}v_k(\omega)
&=0,
\\
\label{ssymm}
s_k(\omega)&=s_{-k}(\tau_k\omega),&
\sum_{k\in\cE}s_k(\omega)&=:s(\omega).
\end{align}
Equation \eqref{divfree} means that $V_k:\Z^d\to [-s^*,s^*]$ is $\pi$-almost surely a bounded and \emph{sourceless flow} on $\Z^d$, or, equivalently, $\Phi:\Z^d\to\R^d$ is a bounded \emph{divergence-free} vector field on $\Z^d$. On the other hand, \eqref{ssymm} implies that
\begin{align}
\label{psiisgrad}
\psi_i(\omega)
=
s_{e_i}(\omega)-s_{e_i}(\tau_{-e_i}\omega),
\qquad
\Psi_i(\omega,x)
=
S_{e_i}(\omega,x)-S_{e_i}(\omega,x-e_i).
\end{align}
That is, the vector field $\Psi:\Z^d\to\R^d$ is component-wise a directional derivative. It follows in particular that
\begin{align}
\label{nodrift1}
\expect{\Psi}=0.
\end{align}
We assume that a similar condition holds for the drift field $\Phi$, too:
\begin{align}
\label{nodrift2}
\expect{\Phi}=\sum_{k\in\cE} k \int_\Omega v_k(\omega) \,{\d}\pi(\omega)=0,
\end{align}
which due to \eqref{divfree}, in the nearest neighbour set-up,  is  obviously the same as assuming that for $k\in\cE$
\begin{align}
\label{nodrift}
\int_\Omega v_k(\omega) \,{\d}\pi(\omega)
=0.
\end{align}
From \eqref{nodrift1} and \eqref{nodrift2} it follows that in the \emph{annealed} mean drift of the walk is nil:
\begin{align}
\notag
\expect{X(t)}
=
\int_\Omega \expectom{X(t)} d\pi(\omega)
=
0.
\end{align}
Under these conditions the law of large numbers
\begin{align}
\label{lln}
\lim_{t\to\infty} t^{-1} X(t) =0,
\qquad
\mathrm{a.s.}
\end{align}
follows directly from the ergodic theorem.

Our next assumption is an ellipticity condition for the symmetric part of the jump rates: there exists another constant  $s_*\in(0,s^*]$ such that for $\pi$-almost all $\omega\in\Omega$ and all $k\in\cE$
\begin{align}
\label{ellipt}
s_k(\omega)\ge s_*, 
\quad 
\pi\text{-a.s.}
\end{align}
Note, that no ellipticity condition is imposed on the jump probabilities $\left(p_k\right)_{k\in\cE}$: it may happen that $p_k=0$ with positive $\pi$-probability. Using a time change we may assume $s_*=1$, and we will occasionally make this assumption for simplicity.

Regarding fluctuations around the law of large numbers \eqref{lln}, we will soon prove that under the ellipticity condition \eqref{ellipt} a \emph{diffusive lower bound} holds: for any fixed vector $v\in\R^d$
\begin{align}
\label{diffusive lower bound}
\varliminf_{t\to \infty} t^{-1} \expect{(v\cdot X(t))^2} >0.
\end{align}
Explicit lower bound will be provided in Proposition \ref{prop:diffusive_bounds} below.

A \emph{diffusive upper bound} also holds under a subtle condition on the covariances of the drift field $\Phi: \Z^d \to \R^d$. Denote
\begin{align}
\label{covariance matrix of Phi}
&
C_{ij}(x)
:=
\cov{\Phi_i(0)}{\Phi_j(x)}
=
\int_\Omega \varphi_i(\omega)\varphi_j(\tau_x\omega)d\pi(\omega),
&&
x\in\Z^d,
\\
\notag
&
D_{ij}(x)
:=
\cov{\Psi_i(0)}{\Psi_j(x)}
=
\int_\Omega \psi_i(\omega)\psi_j(\tau_x\omega)d\pi(\omega),
&&
x\in\Z^d,
\\
\label{Fourier transform of covariance matrix of Phi}
&
\wh C_{ij}(p)
:=
\sum_{x\in\Z^d} e^{\sqrt{-1}x\cdot p} C_{ij}(x),
&&
p\in [-\pi,\pi)^d,
\\
\notag
&
\wh D_{ij}(p)
:=
\sum_{x\in\Z^d} e^{\sqrt{-1}x\cdot p} D_{ij}(x),
&&
p\in [-\pi,\pi)^d.
\end{align}
The Fourier transform is meant as a distribution on
$[-\pi,\pi)^d$. More precisely, by Herglotz's theorem, $\wh C$ and $\wh D$ are positive definite $d\times d$ matrix-valued \emph{measures} on $[-\pi,\pi)^d$. Hence \eqref{nodrift} is equivalent to $\wh C_{ij}(\{0\})=0$, for all $i,j=1,\dots,d$.

The fact that $\Psi$ is a spatial derivative of an $L^2$ function
\eqref{psiisgrad} implies that 
\begin{align}
\label{apriorihcond}
\int_{[-\pi,\pi)^d}
\left(
\sum_{j=1}^d (1-\cos p_j)
\right)^{-1}
\sum_{i=1}^d \wh D_{ii}(p) \, {\d}p <\infty.
\end{align}
A similar \emph{infrared bound} imposed on the covariances of the field $x\mapsto\Phi(x)$ is the  notorious \emph{$\cH_{-1}$-condition} referred to in the title of this paper.

\paragraph{${\cH_{-1}}$-condition}
(first formulation):
We assume
\begin{align}
\label{hcond1}
\int_{[-\pi,\pi)^d}
\left(\sum_{j=1}^d(1-\cos p_j) \right)^{-1} \sum_{i=1}^d \wh C_{ii}(p) \, {\d}p <\infty.
\end{align}

\smallskip
\noindent
For later use we define the positive definite and bounded $d\times d$ matrices
\begin{align}
\label{ctilde}
&
\wt C_{ij}:=
\int_{[-\pi,\pi)^d}
\left(\sum_{j=1}^d (1-\cos p_j) \right)^{-1}
\wh C_{ij}(p)
\, {\d}p <\infty,
\\
\label{dtilde}
&
\wt D_{ij}:=
\int_{[-\pi,\pi)^d}
\left(\sum_{j=1}^d (1-\cos p_j) \right)^{-1}
\wh D_{ij}(p)
\, {\d}p <\infty.
\end{align}

The probabilistic content of the infrared bounds \eqref{apriorihcond} and \eqref{hcond1} is the following. Let $t\mapsto S(t)$ be a continuous time simple symmetric random walk on $\Z^d$ with jump rate $1$, fully independent of the random fields $x\mapsto(\Phi(x), \Psi(x))$. Then \eqref{apriorihcond} and \eqref{hcond1} are in turn equivalent to
\begin{align}
\notag
\lim_{T\to\infty} T^{-1}
\expect{\abs{\int_{0}^{T} \Psi(S(t)){\d}t}^2}
<\infty,
\end{align}
and

\paragraph{${\cH_{-1}}$-condition}
(second formulation):
\begin{align}
\label{hcond2}
&
\lim_{T\to\infty} T^{-1}
\expect{\abs{\int_{0}^{T} \Phi(S(t)){\d}t}^2}
<\infty.
\end{align}

\smallskip
\noindent
The expectations in the last two expressions are  taken over the random walk $t\mapsto S(t)$ \emph{and} the random scenery $x\mapsto(\Phi(x), \Psi(x))$. We omit the straightforward proof of these equivalences. Two more equivalent formulations of the $\cH_{-1}$-condition  \eqref{hcond1}/\eqref{hcond2} will appear later in the paper.

The infrared bounds \eqref{apriorihcond} and \eqref{hcond1} imply a diffusive upper bound: for any fixed vector $v\in\R^d$
\begin{align}
\label{diffusive upper bound}
\varlimsup_{t\to\infty} t^{-1} \expect{(v\cdot X(t))^2} <\infty.
\end{align}
An explicit upper bound will be provided in Proposition \ref{prop:diffusive_bounds} below.

Now, \eqref{diffusive lower bound} and \eqref{diffusive upper bound} jointly \emph{suggest} that the central limit theorem
\begin{align}
\label{clt}
t^{-1/2} X(t) \Rightarrow \cN(0, \sigma^2)
\end{align}
should hold with some non-degenerate $d\times d$ covariance matrix
$\sigma^2$. Attempts to prove the CLT \eqref{clt} under the minimal
conditions of bistochasticity \eqref{bistoch}, ellipticity
\eqref{ellipt}, no drift \eqref{nodrift} and $\cH_{-1}$ \eqref{hcond1} have a notorious history. In Kozlov \cite{kozlov_85} a similar CLT is announced under the somewhat restrictive condition that the random field $x\mapsto P(x)$ in \eqref{jump_probab_field} be finitely dependent. However, as pointed out in Komorowski and Olla \cite{komorowski_olla_03a} the proof in \cite{kozlov_85} is incomplete. In the same paper \cite{komorowski_olla_03a} the CLT \eqref{clt} is stated, but as pointed out in \cite{komorowski_landim_olla_12} this proof is yet again defective. Finally, in \cite{komorowski_landim_olla_12} a complete proof is given, however, with more restrictive conditions: instead of the $\cH_{-1}$-condition \eqref{hcond1} a rather stronger integrability condition on the field $x\mapsto \Phi(x)$  is assumed. See the comments in section \ref{app:Historical remarks}. More detailed historical comments on this story can be found in the notes after chapter 3 of \cite{komorowski_landim_olla_12}. Our main result in the present paper is a complete proof of the CLT \eqref{clt}, under the conditions listed above.

\subsection{Central limit theorem for the random walk}
\label{ss:Central limit theorem for the random walk}

We define the \emph{environment process}, as seen from the random walker:
\begin{align}
\notag
\eta(t):=\tau_{X(t)}\omega
\end{align}
This is a pure jump process on $\Omega$ with bounded total jump rates. So, its construction does not pose any technical difficulty. As already mentioned, it is well known (and easy to check, see e.g.\ Kozlov \cite{kozlov_85}) that due to condition \eqref{bistoch} the probability measure $\pi$ is stationary and ergodic for the Markov process $t\mapsto\eta(t)$. We will denote by $(\cF_t)_{t\ge0}$ the filtration generated by this process:
\begin{align}
\notag
\cF_t:= \sigma(\eta(s): 0\le s\le t).
\end{align}

It is most natural to decompose $X(t)$ as
\begin{align}
\notag
X(t)
=
&
\left\{
X(t)-
\int_{0}^{t} \left(\psi(\eta(s))+\varphi(\eta(s))\right){\d}s
\right\}
+
\int_{0}^{t} \left(\psi(\eta(s))+\varphi(\eta(s))\right){\d}s.
\\[8pt]
\label{martingale decomposition}
=:
&
M(t)+I(t).
\end{align}
In this decomposition the first term is clearly a square integrable
$(\cF_t)$-martingale with stationary and ergodic increments and
conditional covariances (or, quadratic variation)
\begin{align}
\label{qvar}
\condexpect{dM_i(t)dM_j(t)}{\cF_t}=
\delta_{i,j}\left(p_{e_i}(\eta(t))+p_{-e_i}(\eta(t))\right) dt.
\end{align}
Thus, due to the martingale CLT (see e.g.\ \cite{helland_82})
\begin{align}
\notag
t^{-1/2}M(t)\Rightarrow \cN(0,\sigma_M^2),
\end{align}
where
\begin{align}
\label{martvar}
\left(\sigma_M^2\right)_{ij}
=
2 \delta_{i,j} \int_{\Omega}s_{e_i}(\omega) d\pi(\omega).
\end{align}
The difficulty is caused by the compensator integral term $I(t)$.

The following proposition quantifies assertions \eqref{diffusive lower bound} and \eqref{diffusive upper bound}.

\begin{proposition}
\label{prop:diffusive_bounds}
Let $t\mapsto X(t)$ be a random walk in doubly stochastic \eqref{bistoch} random
environment with no drift \eqref{nodrift}. Then the ellipticity \eqref{ellipt} and $\cH_{-1}$ \eqref{hcond1} conditions imply the
following diffusive lower and upper bounds: For any vector
$v\in\R^d$
\begin{align}
\label{diffusive bounds}
2s_* \abs{v}^2
\le
\infsuplim_{t\to\infty} t^{-1} \expect{(v\cdot X(t))^2}
\le
6s^* \abs{v}^2
+
\frac{24}{s_*}\sum_{i,j=1}^d\left(\wt C_{ij}+ \wt D_{ij}\right)v_iv_j,
\end{align}
where $\wt C_{ij}$ and $\wt D_{ij}$ are the matrices defined in \eqref{ctilde} and \eqref{dtilde}.
\end{proposition}

\noindent
The proof of Proposition \ref{prop:diffusive_bounds} is postponed to the next section. Note that the ellipticity condition \eqref{ellipt} is relevant in both (lower and upper) bounds, while the $\cH_{-1}$-condition \eqref{hcond1} is relevant for the upper bound only.

Let us formally state the main result of the present paper.

\begin{theorem}
\label{thm:main}
Let $t\mapsto X(t)$ be a nearest neighbour random walk \eqref{the walk} in random environment, which is bistochastic \eqref{bistoch}, has no drift \eqref{nodrift} and is elliptic \eqref{ellipt}. If in addition the \emph{${\cH_{-1}}$-condition} \eqref{hcond1} holds then

\smallskip
\noindent
(\i)
The asymptotic covariance matrix
\begin{align}
\notag
(\sigma^2)_{ij}
:=
\lim_{t\to\infty} t^{-1}\expect{X_i(t)X_j(t)}
\end{align}
exists, and it is finite and non-degenerate
\begin{align}
\label{bounds on sigmasquared}
2 s_*  I
\le
\sigma^2
\le
6 s^* I_d
+
24 s_*^{-1} \left(\wt C+\wt D\right),
\end{align}
where $I$ is the $d\times d$ unit matrix and $\wt C$, $\wt D$ are the matrices defined in \eqref{ctilde}, \eqref{dtilde}.

\smallskip
\noindent
(\i\i)
Moreover, for any $m\in\N$, $t_1,\dots,t_m\in\R_+$ and any continuous and bounded test function $F:\R^{md}\to\R$
\begin{align}
\notag
\lim_{T\to\infty}
\int_\Omega
\abs{
\expectom{F\left(\frac{X(Tt_1)}{\sqrt{T}}, \dots, \frac{X(Tt_m)}{\sqrt{T}}\right)}
-
\expect{F(W(t_1), \dots, W(t_m))}
} {\d}\pi(\omega)=0,
\end{align}
where $t\mapsto W(t)\in\R^d$ is a Brownian motion with
\begin{align}
\notag
&
\expect{W_i(t)}=0,
&&
\expect{W_i(s)W_j(t)}=\min\{s,t\}(\sigma^2)_{ij}
\end{align}
\end{theorem}

\medskip
\noindent
{\bf Remark} on the jump range of the walk. 
Throughout the paper we speak about nearest neighbour random walk with jump range $\cE$. However, we could consider a more general setup, with jump range $\cU\subset\Z^d$, with the assumptions that (i) $\abs{\cU}<\infty$; (ii) the jump rates are bounded: $p_k(\omega)\le s^*$ almost surely for $k\in\cU$; (iii) the ellipticity condition \eqref{ellipt} holds for a subset $\cU^{\prime}\subset\cU$ which generates $\Z^d$. Under these more general assumptions Theorem 1 remains still valid. The proof remains essentially the same apart of notational changes. 


It is worth noting here that (unlike in the self-adjoint/reversible cases) the $\cH_{-1}$-condition is certainly stronger than assuming just finiteness of the asymptotic variance of the walk, \eqref{diffusive upper bound}. So $\cH_{-1}$ seems to be a sufficient but by no means necessary condition for the CLT to hold. The following question arises very naturally.

\begin{question*}Let $X$ be a stationary, ergodic random walk in a bistochastic random environment, and assume $\expect{|X(t)|^2}\le Ct$. Does it follow that $X$ satisfies a central limit theorem?
\end{question*}

\paragraph*{Structure of the paper.}
The proof of this theorem is the content of sections \ref{s:In the Hilbert space}-\ref{s:The operator B and proof of Theorem 1}. Section \ref{s:In the Hilbert space} contains Hilbert space generalities and most of the notation. Section \ref{s:Relaxed sector condition} describes and slightly extends the relaxed sector condition of \cite{horvath_toth_veto_12} on which we rely. Proofs of the extensions are given in an Appendix (the proofs are similar to those of \cite{horvath_toth_veto_12}, but the statements are stronger). In section \ref{s:The operator B and proof of Theorem 1} we check the conditions of the relaxed sector condition for the concrete case. Remarks, comments (historical and other) and concrete examples are postponed to sections 5-7. 

Let us remark that assuming that $s_k$ is constant for all $k\in\cE$, in other words that the walk is \emph{divergence-free}, removes a number of technical difficulties in the proof. Readers who prefer to see the easier version can see it in the first arxiv version of this paper \cite{circular}

\section{In the Hilbert space \texorpdfstring{$\cL^2(\Omega,\pi)$}{}}
\label{s:In the Hilbert space}

\subsection{Spaces and operators}
\label{ss:Spaces and operators}

It is most natural to put ourselves into the Hilbert space over $\C$
\begin{align}
\notag
\cH:=\Big\{f\in\cL^2(\Omega,\pi):\int_\Omega f {\d}\pi=0\Big\}.
\end{align}
We denote by $T_x$, $x\in\Z^d$, the spatial shift operators
\begin{align}
\notag
T_x f(\omega)
:=
f(\tau_x\omega),
\end{align}
and note that they are unitary:
\begin{align}
\label{shiftops}
T_x^*=T_{-x}=T_x^{-1}.
\end{align}
The $\cL^2$-gradients $\nabla_k$, $k\in\cE$, respectively, $\cL^2$-Laplacian $\Delta$, are:
\begin{align}
\notag
&
\nabla_k
:=
T_k-I,
&&
\nabla_k^*
=
\nabla_{-k},
&&
\norm{\nabla_k}\le 2,
\\[8pt]
\label{Deltaisselfadjoint}
&
\Delta
:=
\sum_{l\in\cE}\nabla_l
=
-
\frac{1}{2}
\sum_{l\in\cE}\nabla_l\nabla_{-l},
&&
\Delta^*
=
\Delta\le0,
&&
\norm{\Delta}\le 4d.
\end{align}
We remark that the norm inequalities above are in fact equalities in
any non-degenerate case, but we will not need this fact.

Due to ergodicity of the $\Z^d$-action $(\Omega, \cF, \pi, \tau_z:z\in\Z^d)$,
\begin{align}
\label{kerDelta}
\Ker(\Delta)=\{0\}.
\end{align}
Indeed, $\Delta f=0$ implies that $0=\sprod{f}{\Delta f} = -\frac12\sum_{k\in\cE} \sprod{\nabla_k f}{\nabla_k f}$ and since all terms are non-negative, they must all be 0 and $f$ must be invariant to translations. Ergodicity to $\Z^d$ actions means that $f$ is constant, and since our Hilbert space is that of functions averaging to zero, $f$ must be zero.

We define the bounded positive operator $\sqd$ in terms of the spectral theorem (applied to the bounded positive operator $\abs{\Delta}:=-\Delta$). Note that due to \eqref{kerDelta} $\Ran \abs{\Delta}$ is dense in $\cH$, and hence so is $\Ran\abs{\Delta}^{1/2}$ which is a superset of it. Hence it follows that $\nsqd :=\left(\sqd\right)^{-1}$ is an (unbounded) positive self-adjoint operator with $\Dom \nsqd = \Ran \sqd$ and $\Ran \nsqd = \Dom \sqd = \cH$. Note that the dense subspace $\Dom \nsqd = \Ran \sqd$ is invariant under, and the operators $\sqd$ and $\nsqd$ commute with the  translations $T_x$, $x\in\Z^d$.

We define the \emph{Riesz operators}: for all $k\in\cE$
\begin{align}
\label{Gamma_k}
\Gamma_k:\Dom \nsqd \to \cH,
\qquad
\Gamma_k=|\Delta|^{-1/2}\nabla_k=\nabla_k|\Delta|^{-1/2},
\end{align}
and note that for any $f\in\Dom \nsqd $
\begin{align}
\notag
\norm{\Gamma_k f}^2
=
\sprod{\nsqd f}{\nabla_{-k}\nabla_k \nsqd f}
\le
\sprod{\nsqd f}{\abs {\Delta} \nsqd f}
=
\norm{f}^2.
\end{align}
Thus, the operators $\Gamma_k$, $k\in\cE$, extend as bounded operators to the whole space $\cH$. The following properties are easy to check:
\begin{align}
\label{Gammaadj}
&
\Gamma_k^*=\Gamma_{-k},
&&
\norm{\Gamma_k}\le 1,
&&
\frac12
\sum_{l\in\cE}\Gamma_l\Gamma_l^*=I.
\end{align}
As before, in fact $\norm{\Gamma_k}=1$ in any non-degenerate case, but
we will not need this fact.

A third equivalent formulation of the ${\cH_{-1}}$-condition \eqref{hcond1}/\eqref{hcond2} is the following:

\paragraph{${\cH_{-1}}$-condition}
(third formulation):
\begin{align}
\label{hcond3}
\varphi_i\in\Dom \nsqd = \Ran \sqd, 
\qquad\qquad
i=1,\dots,d. 
\end{align}
In the case of nearest neighbour walks this is further equivalent to 
\begin{align}
\label{hcond3bis}
v_k\in\Dom \nsqd = \Ran \sqd, 
\qquad\qquad
k\in\cE.
\end{align}

\begin{lemma}
(i) 
Conditions \eqref{hcond3}  and \eqref{hcond2} are equivalent.

\smallskip
\noindent
(ii)
Furthermore, in the case of nearest neighbour walks conditions \eqref{hcond3} and \eqref{hcond3bis} are also equivalent.  
\end{lemma}

\begin{proof}
(i)
Recall that \eqref{hcond2} is formulated in terms of continuous time simple random walk $S$. In operator theory language
\begin{align}
\label{eq:SDelta}
\expectom{\Phi_i(S(t))}
=
e^{t\Delta}\varphi_i(\omega).
\end{align}
Hence
\begin{align*}
\frac{1}{t}\expect{ \abs{\int_{0}^{t} \Phi(S(s)){\,\d}s}^2}
&\stackrel{(*)}{=}
\;\sum_{i=1}^d\int_{0}^{t}
\frac{t-s}{t}
\expect{
2\Phi_i(0)
\Phi_i(S(s))
}
{\,\d}s
\\
&\stackrel{\clap{$\scriptstyle\textrm{(\ref{eq:SDelta})}$}}
{=}
\;\sum_{i=1}^d 2
\int_{0}^{t}
\frac{t-s}{t}
\langle\varphi_i,e^{s\Delta}\varphi_i\rangle
{\,\d}s
\end{align*}
where $(*)$ follows from space stationarity of $\Phi$ (recall that $S$ is independent of $\Phi$, so $\Phi_i(S)$ is just some average of some fixed translations of $\Phi_i$).  An application of the spectral theorem for $|\Delta|$ shows that this is bounded in $t$ if and only if all $\varphi_i\in\Dom|\Delta|^{-1/2}$, $i=1,\dots,d$. 

\smallskip
\noindent 
(ii) 
To conclude from $\varphi\in\Dom|\Delta|^{-1/2}$ that $v\in\Dom|\Delta|^{-1/2}$ we recall that $\varphi_i = (I+T_{-e_i})v_{e_i}=(2I+\nabla_{-e_i})v_{e_i}$. Since $\Gamma_{-e_i}=|\Delta|^{-1/2}\nabla_{-e_i}$ is bounded,  we get that $\nabla_{-e_i} v_{e_i} \in \Dom(|\Delta|^{-1/2})$. Rearranging gives
\begin{align}
\notag
\varphi_i-2v_{e_i} \in \Dom(|\Delta|^{-1/2})
\end{align}
which shows that $\varphi_i\in\Dom(|\Delta|^{-1/2})$ if and only if so is $v_{e_i}$.
\end{proof}

\medskip
\noindent
{\bf Remark.}
Note that equivalence of \eqref{hcond3} and \eqref{hcond3bis} holds \emph{only} in the case of nearest neighbour jumps. If a larger jump range $\cU$ is allowed (see the remark after the formulation of Theorem 1) then \eqref{hcond3bis} is stronger than \eqref{hcond3}. However, the formulation \eqref{hcond3bis} will not be used in the proof of our main result. It will have a role only in the complementary section \ref{app:The stream tensor field}, which is not part of the proof. That part could also be reformulated in the context of finite jump rate, relying only on \eqref{hcond3} but as the main result does not rely on it we will not bother to do that.

Finally, we also define the multiplication operators $M_k, N_k$, $k\in\cE$,
\begin{align}
\label{vmultipl}
&
M_k f(\omega):= v_k(\omega) f(\omega),
&&
M_k^*=M_k,
&&
\norm{M_k}\le s^*,
\\[8pt]
\label{smultipl}
&
N_k f(\omega):= (s_k(\omega)-s_*) f(\omega),
&&
N_k^*=N_k\ge 0,
&&
\norm{N_k}\le s^*
\end{align}
(recall that $s^*$ is the overall \emph{upper bound} on $p$ and $s_*$ is the \emph{lower bound} on the symmetric parts $s$ in the
ellipticity condition (\ref{ellipt})). It is easy to check that the following commutation relations hold due to  \eqref{divfree} and \eqref{ssymm}
\begin{align}
\label{Mnablacommute}
&
\sum_{l\in\cE} M_l \nabla_l
=
-
\sum_{l\in\cE} \nabla_{-l} M_{l},
\\
\notag
&
\sum_{l\in\cE} N_l \nabla_l
=
\sum_{l\in\cE} \nabla_{-l} N_{l}
=
-
\frac12
\sum_{l\in\cE} \nabla_{-l} N_{l} \nabla_{l}.
\end{align}

The \emph{infinitesimal generator} of the stationary environment process $t\mapsto\eta(t)$, acting on the Hilbert space $\cL^2(\Omega,\pi)$ is:
\begin{align}
\notag
Lf(\omega)
= p_k(\omega)(f(\tau_k\omega)-f(\omega)),
\end{align}
which in terms of the operators introduced above is written as
\begin{align}
\label{infgen2}
L=
-D
-T
+A,
\end{align}
with
\begin{align}
\label{opD}
D
&
:=
-s_*\Delta,
\\
\notag
T
&
:=
-
\sum_{l\in\cE} N_l \nabla_l
=
\frac12
\sum_{l\in\cE} \nabla_{-l} N_{l} \nabla_{l},
\\
\label{opA}
A
&
:=
\sum_{l\in\cE} M_l \nabla_l
=
-
\sum_{l\in\cE} \nabla_{-l} M_{l}.
\end{align}
Note that $D=D^*$, $T=T^*$, $A=-A^*$ and
\begin{align}
\label{DboundsT}
0
\le
T
\le
d s^*s_*^{-1}
D.
\end{align}
The inequalities are meant in operator sense. The last one follows from
\begin{align}
\notag
D^{-1/2}TD^{-1/2}
=
\frac{1}{2s_*}\sum_{l\in\cE}\Gamma_{-l}N_l\Gamma_{l},
\end{align}
and hence, due to \eqref{Gammaadj} and \eqref{smultipl}
\begin{align}
\notag
\norm{D^{-1/2}TD^{-1/2}}\le \frac{d s^*}{s_*}
\end{align}
follows, which implies the upper bound in  \eqref{DboundsT}.

\subsection{Proof of Proposition \ref{prop:diffusive_bounds}}
\label{ss:Proof of Proposition 1}

\begin{proof}[Proof of the lower bound in \eqref{diffusive bounds}]
We decompose the displacement process $t\to X(t)$ in such a way that
the forward-and-backward martingale part will be uncorrelated with the rest. The variance of this forward-and-backward martingale will serve as lower bound for the variance of the displacement. Let
\begin{align}
\notag
&
u_k(\omega)
:=
\sgn( v_k(\omega) ) \min\{ \abs{v_k(\omega)} , s_* \},
&&
w_k(\omega)
:=
\sgn( v_k(\omega) ) {\left( \abs{v_k(\omega)} - s_* \right) }_+,
\\[8pt]
\notag
&
q_k(\omega)
:=
s_*+u_k(\omega),
&&
r_k(\omega)
:=
\left( s_k(\omega) - s_* \right) + w_k(\omega).
\end{align}
Note that the skew symmetry \eqref{divfree} of $v_k(\omega)$ is inherited by $u_k(\omega)$ and $w_k(\omega)$:
\begin{align}
\label{vecvec}
u_k(\omega)+ u_{-k}(\tau_k\omega)=0,
\qquad
w_k(\omega)+ w_{-k}(\tau_k\omega)=0.
\end{align}
Further on,
\begin{align}
\notag
&
u_k(\omega)+ w_k(\omega)=v_k(\omega),
&&
q_k(\omega)+ r_k(\omega)=p_k(\omega),
&&
q_k(\omega)\ge0,
&&
r_k(\omega)\ge0.
\end{align}
We further define
\begin{align}
\notag
&
q(\omega):=\sum_{l\in\cE}q_l(\omega)\ge0,
&&
\wt\varphi(\omega):=\sum_{l\in\cE}lq_l(\omega)\in\R^d,
\\
\notag
&
r(\omega):=\sum_{l\in\cE}r_l(\omega)\ge0,
&&
\wt\psi(\omega):=\sum_{l\in\cE}lr_l(\omega)\in\R^d,
\end{align}
and note that
\begin{align}
\notag
&
q(\omega)+r(\omega)=p(\omega),
&&
\wt\varphi(\omega)+\wt\psi(\omega)
=
\varphi(\omega)+\psi(\omega).
\end{align}
Now let $0=\theta_0<\theta_1<\theta_2<\dots$ be the successive jump times of the environment process $t\mapsto \eta(t)$ (or, what is the same, of the random walk $t\mapsto X(t)$):
\begin{align}
\notag
\theta_0:=0,
\qquad
\theta_{n+1}:=\inf\{t>\theta_n: \eta(t)\not=\eta(\theta_n)\},
\end{align}
and define \emph{extra} random variables $\alpha_n\in\{0,1\}$, $n=0, 1,2,\dots$ with the following conditional law, given the trajectory $t\mapsto\eta(t)$: for $N\in\N$ and $a_n\in\{0,1\}$, $n=0,1,\dots,N$,
\begin{align}
\notag
\condprobab{\alpha_n=a_n, \ \ n=0,1,\dots,N}{\eta(t)_{t\ge0}}
=
\prod_{n=0}^N \left(\frac{q(\eta(\theta_{n}))}{p(\eta(\theta_{n}))}\right)^{a_n} \left(\frac{r(\eta(\theta_{n}))}{p(\eta(\theta_{n}))}\right)^{1-a_n}.
\end{align}
In plain words, conditionally on the trajectory $t\mapsto\eta(t)$, the random variables $\alpha_n$, $n=0,1,2,\dots$, are independent biased coin tosses, with probability of head or tail (1 or 0 respectively) equal to the value of $\frac{q(\eta(t))}{p(\eta(t))}$, respectively, $\frac{r(\eta(t))}{p(\eta(t))}$, in the interval $t\in[\theta_n,\theta_{n+1})$. Now, extend piecewise continuously
\begin{align}
\notag
\alpha(t):=\sum_{n=0}^\infty \alpha_n \ind{t\in(\theta_n,\theta_{n+1}]}.
\end{align}
Mind, that $t\mapsto\alpha(t)$ is defined as a \emph{caglad}, not a cadlag process. We decompose the displacement $t\mapsto X(t)$ as follows:
\begin{align}
\notag
X(t)
=
K(t)+L(t)+J(t),
\end{align}
where
\begin{align}
\notag
&
K(t)
:=
\int_0^t \alpha(s) dX(s) - \int_0^t \wt\varphi(s)ds,
\\
\notag
&
L(t)
:=
\int_0^t (1-\alpha(s)) dX(s) - \int_0^t \wt\psi(s)ds,
\\
&
\notag
J(t)
:=\int_0^t \left(\wt\varphi(s) + \wt\psi(s) \right)ds.
\end{align}
Note the following three facts. 

\smallskip
\noindent
(1) \ 
$t\mapsto K(t)$ and $t\mapsto L(t)$, being driven by conditionally independent Poisson flows, are \emph{uncorrelated martingales}, with respect to their own joint filtration. 

\smallskip
\noindent
(2) \ 
$t\mapsto K(t)$  is forward-and-backward martingale with respect to its own past, respectively, future filtration. This is due to  \eqref{vecvec} and to the fact that the symmetric part of its jump rates is constant, $s_*$. Indeed, 
\begin{align*}
\condexpect{K(t+dt)-K(t)}{\eta_t=\omega}
& =
\sum_{l\in\cE} l q_l(\omega) - \wt\varphi(\omega) dt
=
\sum_{l\in\cE} l u_l(\omega) - \wt\varphi(\omega) dt
=
0 dt.
\\
\condexpect{K(t)-K(t-dt)}{\eta_t=\omega}
& =
-\!
\sum_{l\in\cE} l q_{-l}(\tau_l\omega) - \wt\varphi(\omega)dt
=
-\!
\sum_{l\in\cE} l u_{-l}(\tau_l\omega) - \wt\varphi(\omega)dt\\
&=
\sum_{l\in\cE} l u_l(\omega) -\wt\varphi(\omega) dt 
=
0 dt,
\end{align*} 
and hence the claim. 

\smallskip
\noindent
(3) \ 
$t\mapsto J(t)$, being an integral,  is forward-and-backward  predictable with respect to the same filtrations. 

\medskip
\noindent
From these three facts it follows that the process $t\mapsto K(t)$ is \emph{uncorrelated} with $t\mapsto L(t)+J(t)$. Hence, for any vector $v\in\R$
\begin{equation*}
\expect{(v\cdot X(t))^2} =
\expect{(v\cdot K(t))^2}+\expect{(v\cdot (L(t)+J(t)))^2}
\ge
\expect{(v\cdot K(t))^2}
=
2 s_* \abs{v}^2.
\qedhere
\end{equation*}
\end{proof}

\begin{proof}[Proof of the upper bound in \eqref{diffusive bounds}]
We provide upper bounds on the variance of the various terms on the
right hand side of the decomposition $X=M+I$ \eqref{martingale decomposition}.

As shown in \eqref{qvar}-\eqref{martvar} the variance of the martingale term $M(t)$ on the right hand side of \eqref{martingale decomposition} is computed explicitly: for $v\in\R^d$,
\begin{align}
\label{martupp}
\frac{1}{t}\expect{(v\cdot M(t))^2}
=
\sum_{i=1}^d v_i^2 \int_{\Omega}(p_{e_i}(\omega)+p_{-e_i}(\omega)) d\pi(\omega)
\le
2s^*\abs{v}^2.
\end{align}

In order to bound the variance of the integral term $I(t)$ on the
right hand side of \eqref{martingale decomposition} we quote
Proposition 2.1.1 in Olla \cite{olla_01} (alternatively, Lemma
2.4 in \cite{komorowski_landim_olla_12} contains the same result with a different constant). 

\begin{lemma}
\label{lemma:hbound}
Let $t\mapsto\eta(t)$ be a stationary and ergodic Markov process on the probability space $(\Omega,\pi)$, whose infinitesimal generator acting on $\cL^2(\Omega,\pi)$ is $L$. Let $g\in\cL^2(\Omega,\pi)$ such that $\int_\Omega g d\pi=0$. Then
\begin{align}
\notag
\varlimsup_{t\to\infty} \frac{1}{t}
\expect{\max_{0\le s \le t}\abs{\int_0^s g(\eta(u)) {\d}u}^2}
\le
16  \lim_{\lambda\to0}(g, (\lambda I -L-L^*)^{-1}g).
\end{align}
\end{lemma}

\noindent
(Olla denotes the right-hand side by $||g||_{-1}$ --- his definition of
$||g||_{-1}$, (2.1.2) ibid., is different but it is easy to see that
it is equivalent to the above, up to a factor of 2).

The decomposition \eqref{infgen2} of the infinitesimal generator gives
that $-L-L^*\ge 2s_*|\Delta|$, and hence by L\"owner's theorem
(see \cite[Theorem 2.6]{carlen_10} or \cite{lowner_34}) $(-L-L^*)^{-1}\le
1/(2s_*)|\Delta|^{-1}$. It then follows that for any vector $v\in\R^d$
\begin{align}
\label{phiupp}
&
\varlimsup_{t\to\infty}
t^{-1} \expect{\left(\int_0^t v\cdot\varphi(\eta(s))ds\right)^2}
\le
\frac{8}{s_*} ((v\cdot\varphi),|\Delta|^{-1}(v\cdot\varphi))
=
\frac{8}{s_*} \sum_{i,j=1}^d v_i \wt C_{ij} v_j,
\\
\label{psiupp}
&
\varlimsup_{t\to\infty}
t^{-1} \expect{\left(\int_0^t v\cdot\psi(\eta(s))ds\right)^2}
\le
\frac{8}{s_*} ((v\cdot\psi),|\Delta|^{-1}(v\cdot\psi))
=
\frac{8}{s_*} \sum_{i,j=1}^d v_i \wt D_{ij} v_j.
\end{align}
From \eqref{martingale decomposition}, by applying the Cauchy-Schwarz inequality we readily obtain
\[
\expect{\!\left(v\cdot X(t)\right)^2\!}
\le
3
\expect{\!\left(v\cdot M(t)\right)^2\!}
+
3
\expect{\!\left(\int_0^tv\cdot \varphi(\eta(s) ds\right)^{\!2}}
+
3
\expect{\!\left(\int_0^tv\cdot \psi(\eta(s) ds\right)^{\!2}}.
\]
Finally, the upper bound in \eqref{diffusive bounds} follows from here, due to \eqref{martupp}, \eqref{phiupp} and \eqref{psiupp}.
\end{proof}

\section{Relaxed sector condition}
\label{s:Relaxed sector condition}

In  this section we recall and slightly extend the \emph{relaxed sector condition} from \cite{horvath_toth_veto_12}. This is  a functional analytic condition on the operators $D$, $T$ and $A$ from \eqref{infgen2} which ensures that the efficient martingale approximation \`a la Kipnis-Varadhan of integrals of the type of $I(t)$ in \eqref{martingale decomposition} exists. 

A clarification is due here. The relaxed sector condition (Theorem
RSC1 below), is essentially equivalent to the condition that
the range $L\cH_{-1}$ of the infinitesimal generator $L$ be \emph{dense in the $\cH_{-1}$-topology of $\cL^2(\Omega, \pi)$.} (defined by the symmetric part $S:=(L+L^*)/2$ of the infinitesimal generator). This latter one appears in earlier work (see e.g. Olla \cite{olla_01}). But, to the best of  our knowledge it has never been exploited \emph{directly, without stronger sufficient assumptions}. The  \emph{strong} and \emph{graded sector conditions} of Varadhan \cite{varadhan_95}, respectively of Sethuraman, Varadhan and Yau \cite{sethuraman_varadhan_yau_00}, are stronger sufficient conditions for this to hold, and applicable in various circumstances. Nevertheless, the equivalent formulation in \cite{horvath_toth_veto_12} proved to be a very useful one, applicable in conditions where the graded sector condition does not work. In particular, in the context of the present paper. Let us also stress that the graded sector condition itself gets a very transparent and handy proof through the relaxed sector condition. For more details see \cite{horvath_toth_veto_12}. 

Since in the present case the infinitesimal generator $L=-D-T+A$ and
all operators in the decomposition \eqref{infgen2} are \emph{bounded}
we recall the result of \cite{horvath_toth_veto_12} in a slightly restricted form: we do not
have to  worry now about domains and cores of the various operators
$D$, $T$ or $A$.  This section will be fairly abstract.

\subsection{Kipnis-Varadhan theory}
\label{ss:Kipnis-Varadhan theory}

Let $(\Omega, \cF, \pi)$ be a probability space: the state space of a \emph{stationary and ergodic} pure jump Markov process $t\mapsto\eta(t)$ with bounded jump rates. We put ourselves in the complex Hilbert space $\cL^2(\Omega, \pi)$. Denote the \emph{infinitesimal generator} of the semigroup of the process by $L$. Since the process $\eta(t)$ has bounded jump rates the infinitesimal generator $L$ is a bounded operator.  We denote the \emph{self-adjoint} and \emph{skew-self-adjoint} parts of the generator $L$ by
\begin{align}
\notag
S:=-\frac12(L+L^*)\ge0
\qquad
A:=\frac12(L-L^*).
\end{align}
We assume that $S$ is itself ergodic i.e.
\begin{align}
\notag
\Ker(S)=\{c\mathbbm{1} : c\in\C\},
\end{align}
and restrict ourselves to the subspace of codimension 1, orthogonal to the constant functions:
\begin{align}
\notag
\cH:=\{f\in\cL^2(\Omega,\pi): \sprod{\mathbbm{1}}{f}=0\}.
\end{align}
In the sequel the operators $(\lambda I+ S)^{\pm1/2}$, $\lambda\ge0$, will play an important r\^ole. These are defined in terms of the spectral theorem applied to the self-adjoint and positive operator $S$.  The \emph{unbounded} operator $S^{-1/2}$ is self-adjoint on its domain
\begin{align}
\notag
\Dom(S^{-1/2})
=
\Ran(S^{1/2})
=
\{f\in\cH:
\norm{S^{-1/2}f}^2
:=
\lim_{\lambda\to0}\norm{(\lambda I + S)^{-1/2}f}^2
<
\infty
\}.
\end{align}
Let $f\in\cH$. We ask about CLT/invariance principle for the rescaled process
\begin{align}
\label{rescaledintegral}
Y_N(t):=\frac{1}{\sqrt{N}}\int_0^{Nt} f(\eta(s)) {\d} s
\end{align}
as $N\to\infty$.

We denote by $R_\lambda$ the \emph{resolvent} of the semigroup $s\mapsto e^{sL}$:
\begin{align}
\label{resolvent_def}
R_\lambda
:=
\int_0^\infty e^{-\lambda s} e^{sL} {\d} s
=
\big(\lambda I-L\big)^{-1}, \qquad \lambda>0,
\end{align}
and given $f\in\cH$, we will use the notation
\begin{align}
\notag
u_\lambda:=R_\lambda f.
\end{align}

The following theorem is a direct extension to general non-reversible setup of the Kipnis-Varadhan Theorem \cite{kipnis_varadhan_86}. It  yields the \emph{efficient martingale approximation} of the additive functional \eqref{rescaledintegral}.  See T\'oth \cite{toth_86},  or the surveys \cite{olla_01} and \cite{komorowski_landim_olla_12}.

\begin{theoremkv*}
\label{thm:kv}
With the notation and assumptions as before, if the following two limits hold in (the norm topology of) $\cH$:
\begin{align}
\label{conditionA}
\lim_{\lambda\to0}
\lambda^{1/2} u_\lambda=0,
\hskip3cm
\lim_{\lambda\to0} S^{1/2} u_\lambda=v \in\cH,
\end{align}
then
\begin{align}
\notag
\sigma^2
:=
2\lim_{\lambda\to0}\sprod{u_\lambda}{f}
=
2\norm{v}^2
\in
[0,\infty),
\end{align}
exists, and there also exists a zero mean, $\cL^2$-martingale $M(t)$, adapted to the filtration of the Markov process $\eta(t)$, with stationary and ergodic increments and variance
\begin{align}
\notag
\expect{M(t)^2}=\sigma^2t,
\end{align}
such that for $t\in(0,\infty)$
\begin{align}
\notag
\lim_{N\to\infty} 
\expect{\abs{Y_N(t) - \frac{M(Nt)}{\sqrt{N}}}^2} =0.
\end{align}
\end{theoremkv*}

\begin{corollarykv*}
With the same setup and notation,
for any $m\in\N$, $t_1,\dots,t_m\in\R_+$ and $F:\R^{m}\to\R$ continuous and bounded
\begin{align}
\notag
\lim_{N\to\infty}
\int_\Omega
\abs{
\expectom{F(Y_N(t_1), \dots, Y_N(t_m))}
-
\expect{F(W(t_1), \dots, W(t_m))}
}
{\d}\pi(\omega)=0,
\end{align}
where $t\mapsto W(t)\in\R$ is a 1-dimensional Brownian motion with variance $\expect{W(t)^2}=\sigma^2 t$.
\end{corollarykv*}

\subsection{Relaxed sector condition}
\label{ss:Relaxed sector condition}

Let, for $\lambda>0$,
\begin{align}
\label{Clambda_def}
C_\lambda:=(\lambda I + S)^{-1/2} A (\lambda I + S)^{-1/2}.
\end{align}
These are bounded and skew-self-adjoint.

\begin{theoremrsc1*}
\label{thm:rsc1}
Assume that there exist a dense subspace $\cC\subseteq \cH$ and an operator $C:\cC\to \cH$ which is essentially skew-self-adjoint on the core $\cC$ and such that for any vector $\psi\in \cC$ there exists a sequence $\psi_\lambda\in\cH$ such that
\begin{align}
\label{Clambdalimit}
\lim_{\lambda\to 0}\norm{\psi_\lambda-\psi}=0.
\qquad
\text{and}
\qquad
\lim_{\lambda\to 0}\norm{C_\lambda\psi_\lambda-C\psi}=0.
\end{align}
Then for any $f\in\Dom(S^{-1/2})$ the limits \eqref{conditionA} hold and thus the martingale approximation and CLT of Theorem KV follow.
\end{theoremrsc1*}

\paragraph{Remarks}
1. 
The conditions of Theorem RSC1 can be shown to be equivalent to that the sequence of bounded skew-self-adjoint operators $C_\lambda$ converges in the \emph{strong graph limit} sense to the unbounded skew-self-adjoint operator $C$, see Lemma \ref{lem:strrescvg} (ii) below. For various notions of graph limits of operators over Hilbert or Banach spaces see chapter VIII of \cite{reed_simon_vol1_vol2_75}, especially Theorem VIII.26 ibid.
\\
2. 
Theorem RSC1 is a slightly stronger reformulation of Theorem 1 from \cite{horvath_toth_veto_12} where the condition \eqref{Clambdalimit} was slightly stronger. There it was assumed that for any $\varphi\in\cC$, $\lim_{\lambda\to 0}\norm{C_\lambda\varphi-C\varphi}=0$. It turns out that the weaker and more natural condition  \eqref{Clambdalimit} suffices and this has some importance in our next extension, Theorem RSC2. For sake of completeness we give the proof of this theorem in the Appendix.
\\
2.
The operator $C$ is heuristically some version of $S^{-1/2} A S^{-1/2}$. However, it is not sufficient that a naturally densely defined version of $S^{-1/2} A S^{-1/2}$ is skew-Hermitian.  One must  show that its closure is actually skew-self-adjoint. The conditions of Theorem RSC1 require to be careful with domains and with point-wise convergence as $\lambda\to 0$, as above.

\bigskip

RSC refers to \emph{relaxed sector condition}: indeed, as shown in \cite{horvath_toth_veto_12} this theorem contains the \emph{strong sector condition} of \cite{varadhan_95} and the \emph{graded sector condition} of \cite{sethuraman_varadhan_yau_00} as special cases. See the comments at the beginning of Section \ref{s:Relaxed sector condition} for the precise relation of RSC to other sector conditions.  For comments on history, content and variants of Theorem KV we refer the reader to the  monograph \cite{komorowski_landim_olla_12}. For some direct consequences of Theorem RSC1 see \cite{horvath_toth_veto_12}.

Now, we  slightly extend the validity of Theorem RSC1. Assume that the symmetric part of the infinitesimal generator decomposes as
\begin{align}
\notag
S=D+T,
\end{align}
where $D=D^*$, $T=T^*$ and the ``diagonal'' part $D$ dominates $T$ in the following sense: there exists $c<\infty$ so that
\begin{align}
\label{Ddominates}
0\le T \le cD.
\end{align}
Further,  denote
\begin{align}
\label{Blambda_def}
B_\lambda:=(\lambda I + D)^{-1/2} A (\lambda I + D)^{-1/2}.
\end{align}
The following statement is actually a straightforward consequence of Theorem RSC1.

\begin{theoremrsc2*}
\label{thm:rsc2}
Assume that there exist a dense subspace $\cB\subseteq \cH$ and an operator $B:\cB\to \cH$ which is essentially skew-self-adjoint on the core $\cB$  and such that for any vector $\varphi\in \cB$ there exists a sequence $\varphi_\lambda\in\cH$ such that
\begin{align}
\label{Blambdalimit}
\lim_{\lambda\to 0}\norm{\varphi_\lambda-\varphi}=0.
\qquad
\text{and}
\qquad
\lim_{\lambda\to 0}\norm{B_\lambda\varphi_\lambda-B\varphi}=0.
\end{align}
Then for any $f\in\Dom(D^{-1/2})$ the limits \eqref{conditionA} hold and thus the martingale approximation and CLT of Theorem KV follow.
\end{theoremrsc2*}

The proof of Theorem RSC2 is also postponed to the Appendix.

\section{The operator \texorpdfstring{$B=D^{-1/2}AD^{-1/2}$}{B} and
  proof of Theorem \ref{thm:main}}
\label{s:The operator B and proof of Theorem 1}

We apply Theorem RSC2 to our concrete setup, with the operators $D$
and $A$ defined using \eqref{opD} and (\ref{opA}) respectively. Recall that without loss of generality we have fixed $s_*=1$ (see the remark after the ellipticity condition \eqref{ellipt}). Let
\begin{align}
\notag
\cB
:=
\Dom\nsqd
=
\Ran\sqd,
\end{align}
and recall from \eqref{Gamma_k} and \eqref{vmultipl} the definition of the operators $\Gamma_l$ and $M_l$, $l\in\cE$. Define the \emph{unbounded} operator $B:\cB\to\cH$
\begin{align}
\notag
B
:=
-
\sum_{l\in\cE}\Gamma_{-l} M_l \nsqd.
\end{align}
(The definition of $B$ uses our assumption that $s_*=1$, otherwise
with our definitions of $D$ and $A$ we would have needed a factor of
$1/s_*$ before it). 
First we verify \eqref{Blambdalimit}, i.e.\ that $B_\lambda\to B$ \emph{pointwise} on the core $\cB$, where the bounded operator  $B_\lambda$ is expressed by inserting the explicit form of $D$ and $A$, \eqref{opD}, respectively, \eqref{opA}, into the definition \eqref{Blambda_def} of $B_\lambda$:
\begin{align}
\notag
B_\lambda = -\sum_{l\in\cE}(\lambda I - \Delta)^{-1/2}\nabla_{-l}M_l(\lambda I - \Delta)^{-1/2}.
\end{align}
From the spectral theorem for the commutative $C^*$-algebra generated by the shift operators $T_{e_i}$, $i=1,\dots,d$, (see e.g.\ Theorem 1.1.1 on page 2 of \cite{arveson_76}) we obtain that $\|(\lambda I-\Delta)^{-1/2}\nabla_l\|\le 1$, $\|(\lambda I-\Delta)^{-1/2}\abs{\Delta}^{1/2}\|\le 1$, and, moreover, for any $\varphi\in\cH$
\begin{align}
\notag
&
(\lambda I - \Delta)^{-1/2} \nabla_l \varphi
\to
\Gamma_l \varphi
&&
(\lambda I - \Delta)^{-1/2}  \abs{\Delta}^{1/2} \varphi
\to
\varphi, 
&&
\text{as }\lambda\searrow0.
\end{align}
When $\varphi\in\cB$ we get
$(\lambda I-\Delta)^{-1/2}\varphi\to|\Delta|^{-1/2}\varphi$ which allows to write
\begin{align*}
B_\lambda\varphi&=-\sum_{l\in\cE d}(\lambda
I-\Delta)^{-1/2}\nabla_{-l}M_l(\lambda I-\Delta)^{-1/2}\varphi\\
&=-\sum_{l\in\cE d}(\lambda
I-\Delta)^{-1/2}\nabla_{-l}M_l|\Delta|^{-1/2}\varphi+O(\Vert(\lambda I-\Delta)^{-1/2}\varphi
- |\Delta|^{-1/2}\varphi\Vert).
\end{align*}
Hence \eqref{Blambdalimit} follows readily for any
$\varphi\in\cB$. 

With \eqref{Blambdalimit} established, we need to show that $B$ is essentially skew-self-adjoint on $\cB$. We start with a light lemma.

\begin{lemma}
\label{lem:B^*}
(i)
$B:\cB\to\cH$ is skew-Hermitian,
i.e.\ $\langle\varphi,B\psi\rangle=-\langle B\varphi,\psi\rangle$ for
all $\varphi,\psi\in\cB$.
\\
(ii)
The full domain of $B^*$ is
\begin{align}
\label{domain of B^*}
{\cB}^*
=
\{f\in\cH: \sum_{l\in\cE} M_l \Gamma_l f \in\cB\}
\end{align}
and $B^*$ acts on $\cB^*$ by
\begin{align}
\label{B^*}
B^*
:=
-
\nsqd \sum_{l\in\cE} M_l \Gamma_l.
\end{align}
\end{lemma}

\bigskip
\noindent
{\bf Remark:}
It is of crucial importance here that $\cB^*$ in \eqref{domain of B^*}
is the \emph{full} domain of the adjoint operator $B^*$, i.e.~the subspace
of all $f$ such that the linear functional $g\mapsto \sprod{f}{Bg}$ is bounded on
$\cB$. It will not be enough for our purposes just to show that
$\cB^*$ is some core of definition.

\bigskip

\begin{proof}
(i)
Let $f,g\in\cB$. Then, due to \eqref{Mnablacommute}
\begin{align}
\notag
\sprod{f}{Bg}
&
=
-
\sum_{l\in\cE} \sprod{\nsqd f}{\nabla_{-l} M_l \nsqd g}
\\
\notag
&
\stackrel{\clap{$\scriptstyle\textrm{(\ref{Mnablacommute})}$}}{=}
\phantom{-}
\sum_{l\in\cE} \sprod{\nabla_{-l} M_l \nsqd f}{\nsqd g}
=
-\sprod{Bf}{g},
\end{align}
(ii)
Next,
\begin{align}
\notag
\Dom(B^*)
&=
\Big\{
f\in\cH:
(\exists c(f)<\infty)
(\forall g\in\cB):
\Big|\Big\langle f, \sum_{l\in\cE} \Gamma_{-l}M_l \nsqd g\Big\rangle\Big|
\le
c(f)\norm{g}
\Big\}
\\
\notag
&=
\Big\{
f\in\cH:
(\exists c(f)<\infty)
(\forall g\in\cB):
\Big|\Big\langle\sum_{l\in\cE} M_l \Gamma_{l} f,\nsqd g\Big\rangle\Big|
\le
c(f)\norm{g}
\Big\}
\\
\notag
&
=
\Big\{
f\in\cH:
\sum_{l\in\cE} M_l \Gamma_{l} f \in \cB
\Big\},
\end{align}
as claimed. In the last step we used the fact that $\cB$ is the
\emph{full domain} of the self-adjoint operator $\nsqd$. The action \eqref{B^*} of $B^*$ follows from straightforward manipulations.\label{page:other conditions of RSC2}
\end{proof}

Note that Lemma \ref{lem:B^*} in particular implies that
$\cB\subseteq\cB^*$, that $B^*:\cB^*\to\cH$ is in principle an
extension of $-\overline{B}$ and hence the operator $B:\cB\to\cH$ is
closable as the adjoint of any operator is automatically closed. We actually ought to prove that
\begin{align}
\notag
B^*=-\overline{B}.
\end{align}
We apply von Neumann's criterion (see
e.g.\ Theorem VIII.3 of Reed and Simon \cite{reed_simon_vol1_vol2_75}):
If for \emph{some} $\alpha>0$,
\begin{align}
\label{vonNeumann's criterion}
\Ker(B^* \pm \alpha I)=\{0\}
\end{align}
then $B^*=-\overline{B}$.
For reasons which will become clear very soon we will choose $\alpha=s^*$. (Actually any $\alpha\ge s^*$ would work equally well.)
Thus (\ref{vonNeumann's criterion}) is equivalent to showing that the equations
\begin{align}
\label{vonNeumannplus}
\sum_{l\in\cE}M_l \Gamma_l \mu + s^* \sqd \mu = 0
\\
\label{vonNeumannminus}
\sum_{l\in\cE}M_l \Gamma_l \mu - s^* \sqd \mu = 0
\end{align}
admit only the trivial solution $\mu=0$. We will prove this for \eqref{vonNeumannplus}. The other case is done very similarly.

Note that assuming $\mu\in\cB$ the problem becomes fully
trivial. Indeed, inserting $\mu=\abs{\Delta}^{1/2}\chi$ in
\eqref{vonNeumannplus} and taking inner product with $\chi$ we get
\begin{align}
\notag
\sum_{l\in\cE} \langle M_l \nabla_l \chi,\chi\rangle - s^* \langle\Delta \chi,\chi\rangle =0.
\end{align}
The first term is pure imaginary \eqref{Mnablacommute} while the
second term is real \eqref{Deltaisselfadjoint}, giving
that $\langle \Delta\chi,\chi\rangle=0$ which, due to
\eqref{kerDelta}, admits only the trivial solution $\chi=0$. The point is that $\mu$ is not necessarily in $\cB$ so $|\Delta|^{-1/2}\mu$ is not necessarily well defined as an element of $\cH$. Nevertheless, we are able to define a scalar random field $\Psi:\Omega\times\Z^d \to \R$ of stationary increments (rather than stationary) which can be thought of as the lifting of $|\Delta|^{-1/2}\mu$ to the lattice $\Z^d$.

Let, therefore $\mu$ be a putative solution for \eqref{vonNeumannplus} and define, for each $k\in\cE$,
\begin{align}
\label{udef}
u_k:=\Gamma_k \mu.
\end{align}
These are vector components and they also satisfy the \emph{gradient condition}: $\forall \ k,l\in\cE$
\begin{align}
\label{ugradvec}
&
u_k+T_ku_{-k}=0,
&&
u_k+T_ku_l=u_l+T_lu_k.
\end{align}
Note also that
\begin{align}
\notag
\sum_{l\in\cE} u_l = \sqd \mu.
\end{align}
The eigenvalue equation \eqref{vonNeumannplus} becomes
\begin{align}
\label{vonNeumann3}
\sum_{l\in\cE} v_l u_l + s^* \sum_{l\in\cE}u_l=0.
\end{align}
We lift this equation to $\Z^d$. By defining the lattice vector fields $V, U:\Omega\times\Z^d\to\R^d$ as
\begin{align}
\notag
&
V_k(\omega,x):=v_k(\tau_x\omega),
&&
U_k(\omega,x):=u_k(\tau_x\omega),
\end{align}
we obtain the following lifted version of equation \eqref{vonNeumann3}
\begin{align}
\label{vonNeumann3_lifted}
\sum_{l\in\cE} V_l(\omega,x) U_l(\omega,x) + s^* \sum_{l\in\cE}U_l(\omega,x)=0.
\end{align}
Note that $U$ is the $\Z^d$-gradient of a scalar field $\Psi:\Omega\times\Z^d\to\R$, determined uniquely by
\begin{align}
\label{Psidef}
&
\Psi(\omega,0)=0,
&&
\Psi(\omega, x+k) -\Psi(\omega, x) = U_k(\omega, x).
\end{align}
As promised, the scalar field $\Psi$ has stationary increments (or, in
the language of ergodic theory: it is a cocycle), i.e.
\begin{align}
\label{Psistatincr}
\Psi(\omega, y) - \Psi(\omega, x) = \Psi(\tau_x\omega, y-x) - \Psi(\tau_x\omega, 0).
\end{align}
The equation \eqref{vonNeumann3_lifted} gets the form
\begin{align}
\label{harmonic1}
s^*
\sum_{l\in\cE}(\Psi(\omega,x+l)-\Psi(\omega,x))
+
\sum_{l\in\cE}V_l(\omega,x)(\Psi(\omega,x+l)-\Psi(\omega,x))
=0,
\end{align}
Denote the first term by $\lap\Psi$ and the second by $\grad\Psi$
(these are the usual $\Z^d$ Laplacian and gradient, respectively), so the equation becomes
\begin{align}
\label{harmonic2}
s^*  \lap \Psi +  V\cdot \grad \Psi = 0.
\end{align}
We prove that equation \eqref{harmonic1}/\eqref{harmonic2} admits
$\Psi\equiv 0$ as the only solution  satisfying $\expect{\Psi(x)}=0$ for all $x\in\Z^d$. This will
be done using \emph{an auxiliary} random walk in random environment
which will be denoted by $Y$. We remark that in the specific case
where $X$ is \emph{divergence free} i.e.\ $s\equiv 1$, or in general
when $s$ is constant, we get that $Y$ is the same as $X$, but in
general they differ. 

We define the environment for $Y$ on the same probability space
$\Omega$ as $X$. The transfer rates $p^Y_k$, $k\in\cE$ are given by
\begin{align}
\notag
p^Y_k(\omega)= s^* + v_k(\omega).
\end{align}
In other words, we take from $X$ the anti-symmetric part
$v_k=(p_k-p_{-k})/2$ but replace the symmetric part with the constant
$s^*$. 
%
%
The walk $Y$ is also bistochastic, so 
all results proved so far (in particular, stationarity and ergodicity
of $Y$'s environment process $t\mapsto \eta^Y(t):=\tau_{Y(t)}\omega$, and the diffusive lower and upper bounds for $t\mapsto Y(t)$) are in force.

Note that equation \eqref{harmonic1}/\eqref{harmonic2} means exactly
that for a given $\omega\in\Omega$ fixed (that is: in the quenched setup) the field $\Psi(\omega, \cdot):\Z^d\to\R$ is \emph{harmonic} for the random walk $Y(t)$.  Thus the process
\begin{align}
\label{marti}
t\mapsto R(t):=\Psi(Y(t))
\end{align}
is a martingale (with $R(0)=0$) in the quenched filtration $\sigma\left( \omega,Y(s)_{0\le s\le t} \right)_{t\ge0}$. Hence, $t\mapsto R(t)$ is a martingale in its own filtration $\sigma\left( R(s)_{0\le s\le t} \right)_{t\ge0}$, too. We will soon show that $\expect{R(t)^2}<\infty$. From stationarity and ergodicity of the environment process $t\mapsto \eta_t^Y$ and \eqref{Psistatincr} it follows that the process $t\mapsto R(t)$ has stationary and ergodic increments with respect to the annealed measure $\probab{\cdot}:= \int_\Omega \probabom{\cdot} d\pi(\omega)$.
Indeed, let $F(R(\cdot))$ be an arbitrary bounded and measurable functional of the process $t\mapsto R(t)$, $t\ge0$. Using \eqref{Psistatincr}, a straightforward computation shows that
\begin{align}
\notag
\expectom{F(R(t_0+\cdot)-R(t_0))}
=
\expectom{\mathbf{E}_{\eta(t_0)}\left(F(R(\cdot)\right)},
\end{align}
Hence, by stationarity and ergodicity of the environment process $t\mapsto\eta(t)$, the claim follows.

Thus, the process $t\mapsto R(t)$ is a martingale (with $R(0)=0$) with stationary and ergodic increments, in its own filtration $\sigma\left(R(s)_{0\le s\le t} \right)_{t\ge0}$, with respect to the annealed measure $\probab{\cdot}$.

\begin{lemma}
\label{lem:variance of R}
Let $\mu$ be a solution of equation \eqref{vonNeumannplus}, $\Psi$ the harmonic field constructed in \eqref{Psidef} and $R(t)$ the martingale defined in \eqref{marti}. Then
\begin{align}
\label{martingale variance}
\expect{R(t)^2} = 2s^*\norm{\mu}^2 t.
\end{align}
\end{lemma}

\begin{proof}
Since $t\mapsto R(t)$ is a martingale with stationary increments (with respect to the annealed measure $\probab{\cdot}$), we automatically have $\expect{R(t)^2}=\varrho^2 t$ with some $\varrho\ge0$. We now compute $\varrho$.
\begin{align}
\notag
\varrho^2
&
:=
\lim_{t\to0}
\frac{\expect{R(t)^2}}{t}
\stackrel{1}{=}
\lim_{t\to0} \int_\Omega
\frac{\expectom{\Psi(\omega, Y(t))^2}}{t}
\d\pi(\omega)
\\
\notag
&
\stackrel{2}{=}
\int_\Omega \lim_{t\to0}
\frac{\expectom{\Psi(\omega, Y(t))^2}}{t}
\d\pi(\omega)
\stackrel{3}{=}
\sum_{l\in\cE}\int_{\Omega}
\left(s^* + v_l(\omega)\right) \abs{u_l(\omega)}^2
{\d}\pi(\omega)
\\
\notag
&
\stackrel{4}{=}
s^*
\sum_{l\in\cE}\int_{\Omega}
\abs{u_l(\omega)}^2
{\d}\pi(\omega)
\stackrel{5}{=}
s^*
\sum_{l\in\cE}
\norm{\Gamma_l\mu}^2
\stackrel{6}{=}
2 s^* \norm{\mu}^2.
\end{align}
Step 1 is annealed averaging.
Step 2 is easily justified by dominated convergence.
Step 3 drops out from explicit computation of the conditional variance of one jump.
In step 4  we used that due to \eqref{divfree} and \eqref{ugradvec}
$v_{-l}(\omega)\abs{u_{-l}(\omega)}^2 = -
v_{l}(\tau_{-l}\omega)\abs{u_{l}(\tau_{-l}\omega)}^2$ and translation
invariance of the measure $\pi$ on $\Omega$.
In step 5 we use the definition \eqref{udef} of $u_l$.
Finally, in the last step 6 we used the third
identity of \eqref{Gammaadj}.
\end{proof}

\begin{proposition}
\label{prop:Neumann}
The unique solution of \eqref{vonNeumannplus}/\eqref{vonNeumannminus} is $\mu=0$, and consequently the operator $B$ is essentially skew-self-adjoint on the core $\cB$.
\end{proposition}

\begin{proof}
Let $\mu$ be a solution of the equation \eqref{vonNeumannplus}, $\Psi$ the harmonic field constructed in \eqref{Psidef} and $R(t)$ the martingale defined in \eqref{marti}. From the martingale central limit theorem (see e.g.\ \cite{helland_82}) and \eqref{martingale variance} it follows that
\begin{align}
\label{martingale CLT}
\frac{R(t)}{\sqrt{t}}
\Rightarrow \cN(0,2s^*\norm{\mu}^2),
\qquad
\text{ as } t\to\infty.
\end{align}
On the other hand we are going to prove that
\begin{align}
\label{toprob0}
\frac{R(t)}{\sqrt{t}}
\toprob0,
\qquad
\text{ as }t\to\infty.
\end{align}
Jointly, \eqref{martingale CLT} and \eqref{toprob0} clearly imply $\mu=0$, as claimed in the proposition.

The proof of \eqref{toprob0} will combine
\\
(A)
the (sub)diffusive behaviour of the displacement
\begin{align}
\notag
\varlimsup_{T\to\infty}T^{-1}\expect{Y(T)^2}<\infty,
\end{align}
which follows from the ${\cH_{-1}}$-condition,  see \eqref{diffusive bounds}; and
\\
(B)
the fact that the scalar field $x\mapsto \Psi(x)$ having zero mean and stationary increments, cf. \eqref{Psistatincr}, increases \emph{sublinearly} with $\abs{x}$. The sublinearity is the issue here. Since $\Psi$ has stationary, mean zero increments, due to the individual (pointwise) ergodic theorem, it follows that \emph{in any fixed direction} $\Psi$ increases sublinearly almost surely. However, this does not warrant that $\Psi$ increases sublinearly uniformly in $\Z^d$, $d\ge 2$, which is the difficulty we will now tackle.

Let $\delta>0$ and $K<\infty$. Then
\begin{align}
\label{R large}
\probab{\abs{R(t)}>\delta\sqrt{t}}
\le
\probab{\{\abs{R(t)}>\delta\sqrt{t}\} \cap \{\abs{Y(t)}\le K\sqrt{t}\}}
+
\probab{\abs{Y(t)}> K\sqrt{t}}.
\end{align}
From (sub)diffusivity \eqref{diffusive bounds} and Chebyshev's inequality it follows directly that
\begin{align}
\label{X large}
\lim_{K\to\infty}
\varlimsup_{t\to\infty}
\probab{\abs{Y(t)}> K\sqrt{t}}
=0.
\end{align}
We present two proofs of
\begin{align}
\label{R large and X small}
\lim_{t\to\infty}
\probab{\{\abs{R(t)}>\delta\sqrt{t}\} \cap \{\abs{Y(t)}\le K\sqrt{t}\}}
=0,
\end{align}
with $\delta>0$ and $K<\infty$ fixed.  One with bare hands, valid in $d=2$ only, and another one valid in any  dimension which relies on a heat kernel (upper) bound from Morris and Peres  \cite{morris_peres_05}. 

\begin{proof}
[Proof of \eqref{R large and X small} in $d=2$, with bare hands]
We follow here the approach of \cite{berger_biskup_07} where the argument was applied in a different context. In order to keep it short (as another full proof valid in all dimensions follows) we assume separate ergodicity i.e.\ that $(\Omega, \cF, \pi, \tau_{e_i})$ is ergodic for both $i=1,2$. 

First note that
\begin{align}
\label{bound R by Psi}
\probab{\{\abs{R(t)}>\delta\sqrt{t}\} \cap \{\abs{Y(t)}\le K\sqrt{t}\}}
\le
\probab{ \max_{\abs{x}<K\sqrt{t}} \abs{\Psi(x)}> \delta\sqrt{t}}.
\end{align}
Next, since $\Psi$ is harmonic with respect to the random walk $Y(t)$, it obeys the \emph{maximum principle} (this is true for any random walk, no special property of $Y$ is used here). Thus
\begin{align}
\label{maximum principle}
\max_{\abs{x}_\infty\le L}\abs{\Psi(x)}
=
\max_{\abs{x}_\infty = L}\abs{\Psi(x)},
\end{align}
where $\abs{x}_\infty:=\max\{\abs{x_1},\abs{x_2}\}$.  By spatial stationarity
\begin{align}
\label{shiftit}
\begin{gathered}
\max_{\abs{x_1}\le L} \abs{\Psi(x_1,-L)-\Psi(0,-L)}
\sim
\max_{\abs{x_1}\le L} \abs{\Psi(x_1,0)}
\sim
\max_{\abs{x_1}\le L} \abs{\Psi(x_1,+L)-\Psi(0,+L)},
\\
\max_{\abs{x_2}\le L} \abs{\Psi(-L,x_2)-\Psi(-L,0)}
\sim
\max_{\abs{x_2}\le L} \abs{\Psi(0,x_2)}
\sim
\max_{\abs{x_2}\le L} \abs{\Psi(+L,x_2)-\Psi(+L,0)},
\end{gathered}
\end{align}
where $\sim$ stands for equality in distribution. Now, note that $\Psi(x_1,0)$ and $\Psi(0, x_2)$ are Birkhoff sums:
\begin{align}
\notag
\Psi(x_1,0)
=
\sum_{j=0}^{x_1-1} u_{e_1}(\tau_{je_1}\omega),
\qquad
\Psi(0, x_2)
=
\sum_{j=0}^{x_2-1} u_{e_2}(\tau_{je_2}\omega),
\end{align}
where $u_{e_1}(\omega)$ and $u_{e_2}(\omega)$ are zero mean and square integrable (recall the definition of $u$, \eqref{udef}).
Hence, by the ergodic theorem
\begin{align}
\label{ergthm}
L^{-1}
\max\big\{
\max_{\abs{x_1}\le L} \abs{\Psi(x_1,0)},
\max_{\abs{x_2}\le L} \abs{\Psi(0,x_2)}
\big\}
\to0,
\qquad
\text{a.s.,  as }
L\to\infty.
\end{align}
Putting together \eqref{maximum principle}, \eqref{shiftit} and \eqref{ergthm} we readily obtain, for any $\vareps>0$,
\begin{align}
\label{Psi is sublinear}
\lim_{L\to\infty}
\probab{\max_{\abs{x}_\infty\le L}\abs{\Psi(x)}\ge\vareps L}
=0.
\end{align}
Finally, \eqref{R large and X small} follows by applying \eqref{Psi is sublinear} to the right hand side of \eqref{bound R by Psi}.
\end{proof}

\begin{proof}
[Proof of \eqref{R large and X small} in all $d\ge 2$]
We start with the following uniform upper bound on the (quen\-ched) heat kernel of the walk $Y(t)$.

\begin{proposition}
\label{prop:heat kernel bound}
There exists a constant $C=C(d, s^*)$ (depending only on the dimension $d$ and the upper bound $s^*$ on the jump rates) such that for $\pi$-almost all $\omega\in\Omega$ and all $t>0$
\begin{align}
\label{heat kernel bound}
&
\sup_{x\in\Z^d} \probabom{Y(t)=x}\le C t^{-d/2}, 
&&
\pi\text{-a.s.}
\end{align}
\end{proposition}

\begin{proof}
This bound \eqref{heat kernel bound} follows from Theorem 2 of Morris and Peres \cite{morris_peres_05} through Lemma \ref{lem:Morris-Peres}, below, which states essentially the same bound for discrete-time lazy random walks on $\Z^d$ (recall that a random walk is called lazy if there is a lower bound on the probability of the walker staying put at any given point).

\begin{lemma}
\label{lem:Morris-Peres}
Let $V:\Z^d\to[-1,1]^{\cE}$ be a (deterministically given) field such that for all $k\in\cE$ and $x\in\Z^d$
\begin{align}
\label{lifted conditions}
V_k(x)+V_{-k}(x+k)=0,
\qquad
\sum_{l\in\cE}V_l(x)=0.
\end{align}
Define the discrete-time nearest-neighbour, lazy random walk $n\mapsto Y_n$ on $\Z^d$ with transition probabilities
\begin{align}
\label{lazywalk}
\condprobab{Y_{n+1}=y}{Y_n=x}
=
p_{x,y}
:=
\begin{cases}
\frac12
&\text{ if }
y=x,
\\
\frac{1}{4d}(1+V_k(x))
&\text{ if }
y=x+k, \ \  k\in\cE,
\\
0
&\text{ if }
\abs{y-x}>1.
\end{cases}
\end{align}
Then there exists a constant $C=C(d)$ depending only on dimension such that for any $x,y\in\Z^d$
\begin{align}
\label{discrete heat kernel bound}
\condprobab{Y_n=y}{Y_0=x}\le C n^{-d/2}.
\end{align}
\end{lemma}

\begin{proof}
For $A,B\subset\Z^d$, such that $A\cap B=\emptyset$ let
\begin{align}
\notag
Q(A,B)
:=
\sum_{x\in A, y\in B}
p_{x,y}.
\end{align}
For notational reasons we extend the definition of $V_k(x)$, $k\in\cE$, $x\in\Z^d$, as follows
\begin{align}
\notag
V_z(x):=
\begin{cases}
V_k(x)& \text{if } z=k\in\cE,
\\
0 & \text{if } z\not\in\cE.
\end{cases}
\end{align}
For $S\subset\Z^d$, $\abs{S}<\infty$ let $\partial S:=\{(x,y): x\in S, y\in \Z^d\setminus S, \Vert x-y\Vert=1\}$ and note that by the isoperimetric inequality for $\Z^d$ 
\begin{align}
\label{isoperi}
\abs{\partial S}
\ge C \abs{S}^{(d-1)/d}, 
\end{align} 
with some dimension-dependent constant $C$. 
(This discrete isoperimteric inequality is a simple corollary of the classic isoperimetric inequality in $\R^d$. See also Theorem V3.1 in \cite{chavel_01} for a general discretisation result for isoperimteric inequalities.)

We have 
\begin{align}
\notag
Q(S,S^c)
&
=
\sum_{x\in S, y\in S^c} \frac{1}{4d} (1+V_{y-x}(x))
\\
\notag
&
=
\frac{1}{4d} \abs{\partial S}
+
\frac{1}{4d}
\left(
\sum_{x\in S, y\in\Z^d}
V_{y-x}(x)
-
\sum_{x\in S, y\in S}
V_{y-x}(x)
\right)
\\
\label{cut}
&
=
\frac{1}{4d} \abs{\partial S},
\end{align}
where the last equality follows from
\begin{align}
\notag
\sum_{x\in S, y\in\Z^d}
V_{y-x}(x)
&
=
\sum_{x\in S} \sum_{l\in\cE} V_l(x) =0,
\\
\notag
\sum_{x\in S, y\in S}
V_{y-x}(x)
&
=
\frac12
\sum_{x\in S, y\in S}
\left(
V_{y-x}(x)
+
V_{x-y}(y)
\right)
=0,
\end{align}
both of which are consequences of \eqref{lifted conditions}. Yet another consequence of \eqref{lifted conditions} is that the uniform counting measure on $\Z^d$ is stationary to our walk. Hence the isoperimetric profile $\Phi(r)$ (in the sense of Morris and Peres  \cite{morris_peres_05}) is given by
\begin{align}
\notag
\Phi(r)
:=
\inf_{0<\abs{S}\le r}\frac{Q(S,S^c)}{\abs{S}}.
\end{align}
Theorem 2 of \cite{morris_peres_05} (specified to our setup) states that for any $0<\vareps\le 1$, if
\begin{align}
\label{morris-peres thm2}
n> 1 + 4 \int_{4}^{4/\vareps}\frac{{\d} u}{u\Phi^2(u)}
\end{align}
then, for any $x, y \in \Z_d$
\begin{align}
\notag
\condprobab{X_n=y}{X_0=x} \le \vareps.
\end{align}
From \eqref{cut} and the isoperimetric inequality \eqref{isoperi} we have
\begin{align}
\label{isope}
C_1 r^{-1/d} \le \Phi(r) \le C_2 r^{-1/d},
\end{align}
with the constants $0<C_1<C_2<\infty$ depending only on the dimension.
Finally, from \eqref{morris-peres thm2} and \eqref{isope} we readily get \eqref{discrete heat kernel bound}.
\end{proof}

In order to obtain \eqref{heat kernel bound} from \eqref{discrete heat kernel bound}, note that $Y(t)=Y_{\nu(t)}$ where $Y_n$ is a discrete time lazy random walk defined in \eqref{lifted conditions} and \eqref{lazywalk}, with $V_k(x)=v_k(\tau_x\omega)/s^*$ and $t\mapsto\nu(t)$ is a Poisson birth process with intensity $s^*t$ independent of the discrete time walk $Y_n$. Thus
\begin{align}
\notag
\probabom{Y(t)=x}
&
=
e^{-s^*t/2}\sum_{n=0}^\infty \frac{(s^*t/2)^n}{n!} \probabom{Y_n=x}
\\
\notag
&
\le
e^{-s^*t/2}\left(1+\sum_{n=1}^\infty \frac{(s^*t/2)^n}{n!} Cn^{-d/2}\right)
\\
\notag
&
\le
C(d,s^*) t^{-d/2}
\end{align}
This completes the proof of Proposition \ref{prop:heat kernel bound}.
\end{proof}

\subsubsection*{Remarks.}

\newlength{\tempindent}
\setlength{\tempindent}{\parindent}

\begin{enumerate} [leftmargin=0cm,itemindent=0.7cm,labelwidth=\itemindent,labelsep=0cm,align=left,label=(\arabic*)]
\setlength{\parskip}{0cm}\setlength{\parindent}{\tempindent}

\item
The point in Proposition \ref{prop:heat kernel bound} is that it provides uniform upper bound in any (deterministic) environment which satisfies conditions \eqref{lifted conditions}, and thus allows decoupling of the expectation with respect to the walk and with respect to the environment.

\item
In Lemma \ref{lem:Morris-Peres} the ``amount of laziness'' could be any $\delta\in(0,1)$, with appropriate minor changes in the formulation and proof.

\item
Alternative proofs of Proposition \ref{prop:heat kernel bound} are also valid, using either Nash-Sobolev or Faber-Krahn inequalities, see e.g. Kumagai \cite{kumagai_14}. These alternative proofs -- which we do not present here -- are more analytic in flavour. Their advantage is robustness: these proofs are also valid in continuous-space setting (see section \ref{app:Historical remarks} below).

\end{enumerate}

We now return to the proof of \eqref{R large and X small}.  By Chebyshev's inequality
\begin{align}
\label{cheb}
\probab{\{\abs{R(t)}>\delta\sqrt{t}\} \cap \{\abs{Y(t)}\le K\sqrt{t}\}}
\le
\delta^{-2}t^{-1}\expect{\abs{R(t)}^2\ind{\abs{Y(t)}\le K\sqrt{t}}}
\end{align}
Since the scalar field $\Psi$ has stationary increments,
cf. \eqref{Psistatincr}, and zero mean, we get from the $\cL^2$ ergodic theorem that for $k\in\cE$
\begin{align}
\notag
\lim_{n\to\infty}n^{-2}\expect{\abs{\Psi(n k)}^2} =0,
\end{align}
and, consequently,
\begin{align}
\label{Psi sublinear}
\lim_{|x|\to\infty}\abs{x}^{-2}\expect{\abs{\Psi(x)}^2} =0.
\end{align}
Applying in turn the heat kernel bound \eqref{heat kernel bound} of Proposition \ref{prop:heat kernel bound} and the limit \eqref{Psi sublinear} on the right hand side of \eqref{cheb} we obtain
\begin{align}
\notag
t^{-1}\expect{\abs{R(t)}^2 \ind{\abs{Y(t)}\le K\sqrt{t}}}
\le
Ct^{-d/2-1}\sum_{|x|\le K\sqrt{t}}\expect{|\Psi(x)|^2}
\to
0,
\qquad
\textrm{ as }t\to\infty.
\end{align}
Here the first expectation is both on the random walk $Y(t)$ and on the
field $\omega$, while the second is just on the field $\omega$. The point is that with the help of the uniform heat kernel bound of Proposition \ref{prop:heat kernel bound} we can \emph{decouple} the two expectations.

This concludes the proof of \eqref{R large and X small} in arbitrary dimension.
\end{proof}

\noindent
We conclude the proof of the Proposition \ref{prop:Neumann} by noting that from \eqref{R large}, \eqref{X large} and \eqref{R large and X small} we readily get \eqref{toprob0} which, together with \eqref{martingale CLT} implies indeed that $\mu=0$. So \eqref{vonNeumann's criterion} holds with $\alpha=s^*$. We showed that $\Ker(B^*+ s^* I)=\{0\}$, the proof that $\Ker(B^*-s^* I)=\{0\}$ is done in the same way with $Y$ defined using $-V$ instead of $V$. Thus the operator  $B:\cB\to\cH$ is indeed essentially skew-self-adjoint.
\end{proof}

%
%

\begin{proof}[Proof of Theorem \ref{thm:main}]
Proposition \ref{prop:Neumann} verifies that the operator $B$ is essentially skew-self-ad\-joint. The other conditions of Theorem RSC2 are verified \vpagerefrange{s:The operator B and proof of Theorem 1}{page:other conditions of RSC2}. Thus Theorem RSC2 may be applied and we get that for any $f\in\Dom(|\Delta|^{-1/2})$, the time average $\int_0^Nf(\eta(t))$ may be approximated by a Kipnis-Varadhan martingale. The third formulation of the ${\cH_{-1}}$ condition \eqref{hcond3} gives that $v_k\in\Dom(|\Delta|^{-1/2})$ while it is always true that $s_k\in\Dom(|\Delta|^{-1/2})$, \eqref{apriorihcond}. Applying Theorem RSC2 with $f=v_k+s_k$ for each $k\in\{1,\dotsc,d\}$ gives that the compensator $I$ from the decomposition $X=M+I$ \eqref{martingale decomposition} can be approximated with a Kipnis-Varadhan martingale which, we recall, is a stationary martingale $M'$ which is adapted to the filtration of the environment process $\eta$. Hence $M+M'$ is also a stationary martingale and has a CLT. Proposition \ref{prop:diffusive_bounds} gives the bounds \eqref{bounds on sigmasquared}.
\end{proof}

\section{The stream tensor field}
\label{app:The stream tensor field}

The content of this section is not a part of the proof of our main result, but it is an important part of the story and sheds light on the role and limitations of the ${\cH_{-1}}$-condition in this context. We formulate this section in the context of nearest neighbour jumps and part (ii) of Proposition \ref{prop:helmholtz} (below) as presented here relies on the equivalence of \eqref{hcond3} and \eqref{hcond3bis} which is valid only in the nearest neighbour case. However, we remark that this statement, too,  can be easily reformulated for general finite jump rates, but in this case some modifications in the definition of the lattice stream tensor are due and the formulation becomes less transparent. We omit these not particularly instructive details, noting that it is doable with minimum effort.  

The following proposition establishes the existence of the stream tensor field and is essentially Helmholtz's theorem. It is the $\Z^d$ lattice counterpart of Proposition 11.1 from \cite{komorowski_landim_olla_12}. Recall the definition of the field $V:\Omega\times\Z^d\to[-s^*,s^*]^{\cE}$ from \eqref{vlifting}.

\begin{proposition}
\label{prop:helmholtz}
(i)
There exists an antisymmetric tensor field
${H}:\Omega\times\Z^d\to\R^{\cE\times\cE}$ such that
for all $x\in\Z^d$ we have ${H}_{k,l}(\cdot,x)\in\cH$ and
\begin{align}
\label{Thetasymms}
{H}_{l,k}(\omega,x)
=
{H}_{-k,l}(\omega,x+k)
=
{H}_{k,-l}(\omega,x+l)
=
-{H}_{k,l}(\omega,x),
\end{align}
with stationary increments
\begin{align}
\notag
{H}(\omega,y)-{H}(\omega,x)={H}(\tau_x\omega, y-x)-{H}(\tau_x\omega, 0),
\end{align}
such that
\begin{align}
\label{V=curlTheta}
V_k(\omega,x)=\sum_{l\in\cE} {H}_{k,l}(\omega,x).
\end{align}
The realization of the tensor field ${H}$ is uniquely determined by the ``pinning down'' condition \eqref{Theta0} below.

\smallskip
\noindent
(ii)
The ${\cH_{-1}}$-condition \eqref{hcond1} holds if and only if there exist ${h}_{k,l}\in\cH$, $k,l\in\cE$,  such that
\begin{align}
\label{thetatensor}
{h}_{l,k}
=
T_k {h}_{-k,l}
=
T_l {h}_{k,-l}
=
-
{h}_{k,l}
\end{align}
and
\begin{align}
\label{v=curltheta}
v_k(\omega)=\sum_{l\in\cE}{h}_{k,l}(\omega).
\end{align}
In this case the tensor field ${H}$ can be realized as the stationary  lifting of ${h}$:
\begin{align}
\label{lifting of theta}
{H}_{k,l}(\omega,x)={h}_{k,l}(\tau_x\omega).
\end{align}
\end{proposition}

\begin{proof}
(i)
For $k,l,m\in\cE$ define
\begin{align}
\notag
g_{m;k,l}
:=
\Gamma_m\big(\Gamma_l v_k - \Gamma_k v_l\big),
\end{align}
where $\Gamma_l=|\Delta|^{-1/2}\nabla_l$ are the Riesz operators defined in \eqref{Gamma_k}, and note that for all $k,l,m,n\in\cE$
\begin{align}
\label{gistensor}
&
g_{m;l,k}
=
T_kg_{m;-k,l}
=
T_lg_{m;k,-l}
=
-g_{m;k,l},
\\[8pt]
\label{gisgrad}
&
g_{m;l,k}+T_mg_{n;l,k}
=
g_{n;l,k}+T_ng_{m;l,k},
\\[8pt]
\label{curlg=gradv}
&
\sum_{l\in\cE}g_{m;k,l} = \nabla_m v_k.
\end{align}
\eqref{gistensor} means that that keeping the subscript $m\in\cE$ fixed, $g_{m;k,l}$ has exactly the symmetries of a $\cL^2$-tensor variable indexed by  $k,l\in\cE$. \eqref{gisgrad} means that, on the other hand, keeping $k,l\in\cE$ fixed, $g_{m;k,l}$ is a $\cL^2$-gradient in the subscript $m\in\cE$. Finally, \eqref{curlg=gradv} means that the $\cL^2$-divergence of tensor $g_{m;\cdot,\cdot}$ is actually the $\cL^2$-gradient of the vector $v_\cdot$.

Let  $G_{m;k,l}:\Omega\times\Z^d\to\R$ be the lifting  $G_{m;k,l}(\omega,x):=g_{m;k,l}(\tau_x\omega)$. By \eqref{gisgrad}, for any $k,l\in\cE$ fixed $\left(G_{m;k,l}(\omega,x)\right))_{m\in\cE}$ is a lattice gradient. The increments of ${H}_{k,l}$ are defined by
\begin{align}
\label{gradTheta=G}
{H}_{k,l}(\omega,x+m)-{H}_{k,l}(\omega,x)
=
G_{m;k,l}(\omega,x),
\ \ \ m\in\cE.
\end{align}
This is consistent, due to \eqref{gisgrad}.

Next, in order to uniquely determine the tensor field $H$,  we ``pin down'' its values at $x=0$.
For $e_i,e_j\in\cE_{+}$ choose
\begin{align}
\label{Theta0}
\begin{aligned}
&
{H}_{e_i,e_j}(\omega,0)
=
0,
&&
{H}_{-e_i,e_j}(\omega,0)
=
-g_{-e_i;e_i,e_j}(\omega),
\\[8pt]
&
{H}_{e_i,-e_j}(\omega,0)
=
g_{-e_j;e_i,e_j}(\omega),
&&
{H}_{-e_i,-e_j}(\omega,0)
=
-g_{-e_i;e_i,e_j}(\omega)+g_{-e_j;e_i,e_j}(\tau_{-e_i}\omega).
\end{aligned}
\end{align}
The tensor field ${H}$ is fully determined by \eqref{gradTheta=G} and \eqref{Theta0}. Due to \eqref{gistensor} and \eqref{curlg=gradv}, \eqref{Thetasymms}, respectively, \eqref{V=curlTheta} will hold, indeed.

(ii)
We show equivalence with $v_k\in\Dom(|\Delta|^{-1/2})$ \eqref{hcond3}. 
First we prove the \emph{only} if part. Assume \eqref{hcond3} and let
\begin{align}
\notag
{h}_{k,l}
=
\Gamma_l \abs{\Delta}^{-1/2} v_k - \Gamma_k  \abs{\Delta}^{-1/2} v_l
=
\abs{\Delta}^{-1/2} \big( \Gamma_l  v_k - \Gamma_k  v_l\big).
\end{align}
Hence \eqref{thetatensor} and \eqref{v=curltheta} are readily obtained. Next we prove the if part. Assume that there exist $h_{k,l}\in\cH$ with the symmetries \eqref{thetatensor} and $v_k$ is realized as in \eqref{v=curltheta}. Then we have
\begin{align}
\notag
v_k
=
\sum_{l\in\cE} h_{k,l}
=
\frac12 \sum_{l\in\cE} \left(h_{k,l} + h_{k,-l}\right)
=
-\frac12 \sum_{l\in\cE} \nabla_l h_{k,-l}
=
-\frac12 \abs{\Delta}^{1/2} \sum_{l\in\cE} \Gamma_l h_{k,-l},
\end{align}
which shows indeed \eqref{hcond3}.
\end{proof}

\paragraph{${\cH_{-1}}$-condition}
(fourth formulation):
The drift vector field $V$ is realized as the curl of a \emph{stationary and square integrable}, zero mean  tensor field ${H}$, as shown in \eqref{V=curlTheta}.

\bigskip
\noindent
{\bf Remark.}
If the $\cH_{-1}$-condition \eqref{hcond1} does not hold it may still be possible that there exists a \emph{non-square integrable} tensor variable ${h}:\Omega\to\R^{\cE\times\cE}$ which has the symmetries \eqref{thetatensor} and with $v:\Omega\to\R^{\cE}$ realized as in \eqref{v=curltheta}. Then let ${H}:\Omega\times\Z^d\to\R^{\cE\times\cE}$ be the stationary lifting \eqref{lifting of theta} and we still get \eqref{V=curlTheta} with a \emph{stationary but not square integrable} tensor field. Note that this is not decidable in terms of the covariance matrix \eqref{covariance matrix of Phi} or its Fourier transform \eqref{Fourier transform of covariance matrix of Phi}. The question of diffusive (or super-diffusive) asymptotic behaviour of the walk $t\mapsto X(t)$ in these cases is fully open.

\medskip

In the next proposition -- which essentially follows an argument from Kozlov  \cite{kozlov_85} -- we give a sufficient condition for the ${\cH_{-1}}$-condition \eqref{hcond1} to hold.

\begin{proposition}
\label{prop:suff cond for H-1}
If $p\mapsto \wh C(p)$ is twice continuously differentiable function in a neighbourhood of $p=0$ then the ${\cH_{-1}}$-condition \eqref{hcond1} holds.
\end{proposition}

\begin{proof}
For the duration of this proof we introduce the notation
\begin{align}
\notag
B_{k,l}(x):=\expect{V_k(0)V_l(x)},
\quad
\wh B_{k,l}(p) := \sum_{x\in \Z^d} e^{\sqrt{-1}x\cdot p} B_{k,l} (x),
\end{align}
with $k,l\in\cE, x\in\Z^d, p\in[-\pi,\pi]^d$. Hence for $i,j\in\{1,\dots,d\}$
\begin{align}
\notag
\wh C_{ij}(p)
=
\wh B_{e_i,e_j}(p) - \wh B_{-e_i,e_j}(p) - \wh B_{e_i,-e_j}(p) + \wh B_{-e_i,-e_j}(p).
\end{align}
(The identity is meant in the sense of distributions.)

Note that  due to the first clause in \eqref{lifted conditions}
\begin{align}
\label{Chatvector}
\wh B_{k,l}(p)
=
-e^{\sqrt{-1}p\cdot k}\wh B_{-k,l}(p)
=
-e^{-\sqrt{-1}p\cdot l}\wh B_{k,-l}(p)
=
e^{\sqrt{-1}p\cdot (k-l)}\wh B_{-k,-l}(p).
\end{align}
Using \eqref{Chatvector} in the above expression of $C(p)$ in terms of $B(p)$, direct computations yield
\begin{align}
\notag
\wh C_{ij} = \left(1+e^{-\sqrt{-1}p\cdot e_i}\right) \left(1+e^{\sqrt{-1}p\cdot e_j}\right) \wh B_{e_i,e_j}(p).
\end{align}
Thus, the regularity condition imposed on $p\mapsto C(p)$ is equivalent to assuming the same regularity about $p\mapsto \wh B(p)$.

Next, due to the second clause of \eqref{lifted conditions}
\begin{align}
\label{Chatdivfree}
\sum_{k\in\cE}\wh B_{k,l}(p)
=
\sum_{l\in\cE}\wh B_{k,l}(p)
=
0,
\end{align}
and, from \eqref{Chatvector} and \eqref{Chatdivfree} again by direct computations we obtain
\begin{align}
\label{equiv0}
\sum_{k,l\in\cE}
(1-e^{-\sqrt{-1}p\cdot k})(1-e^{\sqrt{-1}p\cdot l})\wh B_{k,l}(p)\equiv0.
\end{align}
At $p=0$ we apply $\partial^2/\partial p_i\partial p_j$ to \eqref{equiv0} and get
\begin{align}
\label{C(0)=0}
\wh C_{ij}(0)
=
\sum_{k,l\in\cE}
k_il_j \wh B_{k,l}(0)=0,
\qquad
i,j=1,\dots,d.
\end{align}
Since $\wh C_{j,i}(p) = \wh C_{ij}(-p) = \overline{\wh C_{ij}(p)} $
and $p\mapsto \wh C(p)$ is assumed to be twice continuously differentiable at $p=0$, from \eqref{C(0)=0} it follows that
\begin{align}
\notag
\wh C(p)=\Ordo(\abs{p}^2),
\qquad
\text{ as }
\abs{p}\to0,
\end{align}
which implies \eqref{hcond1}.
\end{proof}

\noindent
In particular it follows that sufficiently fast decay of correlations of the divergence-free drift field $V(x)$ implies the ${\cH_{-1}}$-condition \eqref{hcond1}. Note that the divergence-free condition \eqref{divfree} is crucial in this argument.

\section{Historical remarks}
\label{app:Historical remarks}

There exist a fair number of important earlier results to which we should compare Theorem \ref{thm:main}.

\begin{enumerate} [leftmargin=0cm,itemindent=0.7cm,labelwidth=\itemindent,labelsep=0cm,align=left,label=(\arabic*)]
\setlength{\parskip}{0cm}\setlength{\parindent}{\tempindent}

\item
In Kozlov \cite{kozlov_85}, Theorem II.3.3 claims the same result under the supplementary restrictive condition that the random field of jump probabilities $x\mapsto P(x)$ in \eqref{jump_probab_field} be \emph{finitely dependent}. However, as pointed out by Komorowski and Olla \cite{komorowski_olla_02}, the proof is incomplete there. Also, the condition of finite dependence of the field of jump probabilities is a very serious restriction.

\item
In Komorowski and Olla \cite{komorowski_olla_03a}, Theorem 2.2, essentially the same result is announced as above. However, as noted in section 3.6 of \cite{komorowski_landim_olla_12} this proof is yet again incomplete.

\item
To our knowledge the best fully proved result is Theorem 3.6 of \cite{komorowski_landim_olla_12} where the same result is proved under the condition that the stream tensor field $x\mapsto H(x)$ of Proposition \ref{prop:helmholtz} be stationary and in $\cL^{\max\{2+\delta, d\}}$, $\delta>0$, rather than $\cL^2$. Note that the conditions of our theorem only request that the tensor field $x\mapsto H$ be square integrable. The proof of Theorem 3.6 in \cite{komorowski_landim_olla_12} is very technical, see sections 3.4 and 3.5 of the monograph.

\item
The special case when the tensor field ${H}$ is actually in $\cL^{\infty}$ is fundamentally simpler. In this case the so-called \emph{strong sector condition} of Varadhan \cite{varadhan_95} applies directly. This was noticed in \cite{komorowski_olla_03a}. See also section 3.3 of \cite{komorowski_landim_olla_12} and section \ref{app:Examples} below.

\item
Examine the following \emph{diffusion} problem is as follows. Let $t\mapsto X(t)\in \R^d$ be the strong solution of the SDE
\begin{align}
\label{sde}
{\d}X(t)
=
{\d}B(t) + \Phi(X(t)){\d}t,
\end{align}
where $B(t)$ is standard $d$-dimensional Brownian motion and $\Phi:\R^d\to\R^d$ is a stationary and ergodic (under space-shifts) vector field on $\R^d$ which has zero mean
\begin{align}
\notag
\expect{\Phi(x)}=0,
\end{align}
and is almost surely \emph{divergence-free}:
\begin{align}
\label{divfree in continuous space}
\div \Phi \equiv 0,
\ \ \
\mathrm{a.s.}
\end{align}
It is analogous to the discrete-space problem studied in this paper in the case that $s_k$ is constant for all $k\in\cE$. In this case the $\cH_{-1}$-condition is
\begin{align}
\label{hcond in continuous space}
\sum_{i=1}^d\int_{\R^d} \abs{p}^{-2}  \wh C_{ii}(p) {\d}p <\infty,
\end{align}
where
\begin{align}
\notag
\wh C_{ij}(p)
:=
\int_{\R^d} \expect{\Phi_i(0)\Phi_j(x)}  e^{\sqrt{-1}p\cdot x} {\d}x,
\qquad
p\in\R^d.
\end{align}
It is a fact that, similarly to the $\Z^d$ lattice case, under minimally restrictive regularity conditions, a stationary and square integrable divergence-free drift field $x\mapsto \Phi(x)$ on $\R^d$ can be written as the curl of an antisymmetric stream tensor field with stationary increments $H:\R^d\to\R^{d\times d}$:
\begin{align}
\notag
\Phi_i(x) = \sum_{j=1}^d \frac{\partial H_{ji}}{\partial x_j}(x).
\end{align} 
This is essentially Helmholtz's theorem. 
See Proposition 11.1 of \cite{komorowski_landim_olla_12}, which is the continuous-space analogue of Proposition \ref{prop:helmholtz} of section \ref{app:The stream tensor field} above. As shown in \cite{komorowski_landim_olla_12}, the ${\cH_{-1}}$-condition \eqref{hcond in continuous space} is equivalent with the fact that the stream tensor ${H}$ is \emph{stationary} (not just of stationary
increments) \emph{and square integrable}. The case of bounded ${H}$ was first considered in Papanicolaou and Varadhan  \cite{papanicolaou_varadhan_81}. This paper is historically the first instant where the problem of diffusion in stationary divergence-free drift field was considered with mathematical rigour. Homogenization and central limit theorem for the diffusion \eqref{sde}, \eqref{divfree in continuous space} in \emph{bounded} stream field, ${H}\in\cL^{\infty}$, was first proven in Osada \cite{osada_83}. Today the strongest result in the continuous space-time setup is due to Oelschl\"ager \cite{oelschlager_88} where homogenization and CLT for the displacement is proved for square-integrable stationary stream tensor field, ${H}\in\cL^2$.  Oelschl\"ager's proof consists in truncating the stream tensor and bounding the error.  If the stream tensor field is stationary \emph{Gaussian} then -- as noted by Komorowski and Olla \cite{komorowski_olla_03b} -- the \emph{graded sector condition} of \cite{sethuraman_varadhan_yau_00} can be applied. See also chapters 10 and 11 of \cite{komorowski_landim_olla_12} for all existing results on the diffusion model \eqref{sde}, \eqref{divfree in continuous space}. 

\item
Attempts to apply Oelschl\"ager's method in the discrete ($\Z^d$ rather than $\R^d$) setting run into enormous technical difficulties, see chapter 3 of \cite{komorowski_landim_olla_12} and seemingly this approach can't be fully accomplished beyond the overly restrictive condition ${H}\in\cL^{\max\{2+\delta, d\}}$. The main result of this paper, Theorem \ref{thm:main} fills this gap between the restrictive condition ${H}\in\cL^{\max\{2+\delta, d\}}$ of Theorem 3.6 in \cite{komorowski_landim_olla_12} and the minimal restriction ${H}\in\cL^{2}$. The content of our Theorem \ref{thm:main} is the discrete $\Z^d$-counterpart of Theorem 1 in Oelschl\"ager \cite{oelschlager_88}. We also stress that our proof is conceptually and technically much simpler that of Theorem 3.6 in \cite{komorowski_landim_olla_12} or Theorem 1 in \cite{oelschlager_88}. The continuous space-time diffusion model --- under the same regularity conditions as those of Oelschl\"ager \cite{oelschlager_88} can be treated in a very similar way reproducing this way Theorem 1 of \cite{oelschlager_88} in a conceptually and technically simpler way. In order to keep this paper relatively short and transparent, those details will be presented elsewhere.

\item
There exist results on \emph{super-diffusive} behaviour of the random walk in doubly stochastic random environment \eqref{the walk}, \eqref{bistoch} or diffusion in divergence-free random drift field \eqref{sde}, \eqref{divfree in continuous space}, when the ${\cH_{-1}}$-condition \eqref{hcond1} fails to hold. In Komorowski and Olla \cite{komorowski_olla_02} and  T\'oth and Valk\'o \cite{toth_valko_12} the diffusion model \eqref{sde}, \eqref{divfree in continuous space} is considered when the drift field $\Phi$ is Gaussian and the stream tensor field ${H}$ is \emph{genuinely} delocalized: of stationary increment but not stationary. Super-diffusive bounds are proved.

\end{enumerate}

\section{Examples}
\label{app:Examples}

Before formulating concrete examples let us spend a few words about the physical motivation and phenomenology of the problem considered. The continuous case discussed in the previous section, diffusion in divergence-free drift field, cf.  \eqref{sde}-\eqref{divfree in continuous space} may model the drifting of a suspended particle in stationary turbulent incompressible flow. Very similarly, the lattice counterpart \eqref{the walk} with jump rates satisfying \eqref{bistoch} describe a random walk whose local drift is driven by a stationary source- and sink-free flow. The interest in the asymptotic description of this kind of displacement dates back to the discovery of turbulence. However, divergence-free environments appear in many other natural contexts, too. See e.g.\ \cite[chapter 11]{komorowski_landim_olla_12} or a surprising recent application to group theory by Bartholdi and Erschler \cite{bartholdi_erschler_11}.

A phenomenological picture of these walks can be formulated in terms of randomly oriented cycles. Imagine that a translation invariant random ``soup of cycles'' --- that is, a Poisson point process of oriented cycles --- is placed on the lattice, and the walker is drifted along by these whirls. Now, local small cycles contribute to the diffusive behaviour. But occasionally very large cycles may cause on the long time scale faster-than-diffusive transport. Actually, this happens: in Komorowski and Olla \cite{komorowski_olla_02} and T\'oth and Valk\'o \cite{toth_valko_12} anomalous \emph{superdiffusive} behaviour is proved in particular cases when the $\cH_{-1}$-bound \eqref{hcond in continuous space} doesn't hold. Our result establishes that on the other hand, the $\cH_{-1}$-bound \eqref{hcond1} ensures not only boundedness of the diffusivity but also normal behaviour under diffusive scaling.

And now, to some examples: 

\setlength{\tempindent}{\parindent}
\begin{enumerate} 
[leftmargin=0cm,itemindent=0.7cm,labelwidth=\itemindent,labelsep=0cm,align=left,label=(\arabic*)]
\setlength{\parskip}{0cm}\setlength{\parindent}{\tempindent}

\item
\emph{Stationary and bounded stream field:}
When there exists a \emph{bounded} tensor valued variable ${h}:\Omega\to\R^{\cE\times\cE}$ with the symmetries \eqref{thetatensor} and such that \eqref{v=curltheta}  holds we define the multiplication operators $M_{k,l}$, $k,l\in\cE$, acting on $f\in\cH$:
\begin{align}
\label{thetamultipl}
M_{k,l} f(\omega):= {h}_{k,l}(\omega) f(\omega).
\end{align}
These are bounded selfadjoint  operators and they inherit the symmetries of ${h}$ (recall the shift operators $T_k$, $k\in\cE$ from \eqref{shiftops}):
\begin{align}
\label{Nsymmetries}
\begin{gathered}
M_{l,k}=T_k M_{-k,l}T_{-k}=T_l M_{k,-l}T_{-l}=-M_{k,l},
\\[8pt]
\sum_{l\in\cE}M_{k,l}=M_k.
\end{gathered}
\end{align}
As an alternative to \eqref{opA}, using \eqref{Nsymmetries}, the skew-self-adjoint part of the infinitesimal generator is expressed as
\begin{align}
\label{Aalt}
A=\sum_{k,l\in\cE} \nabla_{-k}M_{k,l}\nabla_l.
\end{align}
In \cite{komorowski_olla_03a} and \cite{komorowski_landim_olla_12} this form of the operator $A$ is used. The operators $M_{k,l}$ are bounded and so is the operator
\begin{align}
\label{bibabu}
B
:=
\abs{\Delta}^{-1/2} A \abs{\Delta}^{-1/2}
=
\sum_{k,l\in\cE} \Gamma_{-k}M_{k,l}\Gamma_l
\end{align}
which plays a key r\^ole in our proof. Due to boundedness of $B$ the \emph{strong sector condition} is valid in these cases and the central limit theorem for the displacement readily follows. See \cite{komorowski_olla_03a} and section 3.3 of \cite{komorowski_landim_olla_12}.

Finitely dependent constructions of this type appear in Kozlov \cite{kozlov_85}. The so-called \emph{cyclic walks} analysed in \cite{komorowski_olla_03a} and in section 3.3 of \cite{komorowski_landim_olla_12} are  also of this nature. 

When the tensor variables ${h}:\Omega\to\R^{\cE\times\cE}$ in \eqref{thetatensor} are in $\cL^2\setminus \cL^\infty$, the multiplication operators $M_{k,l}$ defined in \eqref{thetamultipl} are \emph{unbounded}, the representation \eqref{Aalt} of the skew-self-adjoint part of the infinitesimal generator and the operator $B$ defined in \eqref{bibabu} become just \emph{formal}. Nevertheless,  Theorem 1 in Oelschl\"ager \cite{oelschlager_88} and theorem 3.6 in \cite{komorowski_landim_olla_12} are proved by controlling approximations of $h_{k,l}$ and the unbounded operators $M_{k,l}$ by truncations at high levels.

\item
\emph{Stationary, square integrable but unbounded stream field:}
We let, in arbitrary dimension $d$, $\Psi:\Z^d+(1/2,\dots, 1/2)\to\Z$ be a stationary, scalar, Lipschitz field with Lipschitz constant 1. As shown in Peled  \cite{peled_10}, such fields exist in sufficiently high dimension. Define ${H}:\Z^d\to\R^{\cE_{2}\times\cE_{2}}$ by
\begin{align}
\notag
&
{H}_{e_i,e_j}(x)
:=
\frac 1{d}\Psi(x+(e_i+e_j)/2),
&&
x\in\Z^d,
\ \
1\le i<j \le d,
\end{align}
and extend to $\left({H}_{k,l}(x)\right)_{k,l\in\cE}$ by the symmetries \eqref{Thetasymms}. The tensor field ${H}:\Z^d\to\R^{\cE_{2}\times\cE_{2}}$ defined this way will be stationary and $\cL^2$, but not necessary in $\cL^\infty$ --- the uniform graph homomorphism of Peled \cite{peled_10}, for example, is not bounded. Nevertheless, $V$ is bounded by 1, as it should, since $|H_{k,l}(x)+H_{-k,l}(x)|=|H_{k,l}(x)-H_{k,l}(x-k)|\le \frac 1d$ and $V$ is a sum of $d$ such terms.

\item
\emph{Randomly oriented  Manhattan lattice:}
Let $u_i:\Z^{d-1}\to \{-1,+1\}$, $i=1,\dots, d$,  be translation invariant and ergodic, zero mean random fields, which are independent between them. Denote their covariances
\begin{align}
\notag
c_i(y)
&
:=
\expect{u_i(0), u_i(y)},
&&
y\in\Z^{d-1},
\\[8pt]
\notag
\hat c_i(p)
&
:=
\sum_{y\in\Z^{d-1}} e^{\sqrt{-1}p\cdot y} c_i(y),
&&
p\in[-\pi,\pi)^{d-1}.
\end{align}
Define now the lattice vector field
\begin{align}
\notag
V_{\pm e_i}(x)
:=
\pm  u_i(x_1, \dots, x_{i-1}, \cancel{x_i}, x_{i+1}, \dots, x_d).
\end{align}
Then the random vector field $V$ will satisfy  all conditions in \eqref{lifted conditions} and $t\mapsto X(t)$ will actually be a random walk on the lattice $\Z^d$ whose line-paths parallel to the coordinate axes are randomly oriented in a shift-invariant and ergodic way. This oriented graph is called the \emph{randomly oriented Manhattan lattice}. The covariances $C$ and $\wh C$ defined in \eqref{covariance matrix of Phi}, respectively, \eqref{Fourier transform of covariance matrix of Phi} will be
\begin{align}
\notag
C_{ij} (x)
&
=
\delta_{i,j} c_i(x_1, \dots, x_{i-1}, \cancel{x_i}, x_{i+1}, \dots, x_d),
\\[8pt]
\notag
\wh C_{ij}(p)
&
=
\delta_{i,j} \delta(p_i) \hat c_i(p_1, \dots, p_{i-1}, \cancel{p_i}, p_{i+1}, \dots, p_d).
\end{align}
The $\cH_{-1}$-condition \eqref{hcond1} is in this case
\begin{align}
\label{hcond for manhattan}
\sum_{i=1}^d
\int_{[-\pi,\pi]^{d-1}}
\wh D(q)^{-1}
\hat c_{i}(q) {\d}q <\infty.
\end{align}
In the particular case when the random variables $u_i(y)$, $i\in\{1,\dots, d\}$, $y\in\Z^{d-1}$, are \emph{independent fair coin-tosses}, $\hat c_i(q)\equiv1$. In this case, for $d=2,3$ the ${\cH_{-1}}$-condition  \eqref{hcond for manhattan} fails to hold, the tensor field ${H}$ is \emph{genuinely} of stationary increments. In these cases super-diffu\-si\-vi\-ty of the walk $t\mapsto X(t)$ can be proved with the method of Tarr\`es, T\'oth and Valk\'o \cite{tarres_toth_valko_12} (in the $2d$ case), respectively, of T\'oth and Valk\'o \cite{toth_valko_12} (in the $3d$ case). In dimensions $d\ge4$ the ${\cH_{-1}}$-condition  \eqref{hcond for manhattan} (and thus \eqref{hcond1}) holds and the central limit theorem for the displacement follows from our Theorem \ref{thm:main}.

\end{enumerate}

\section{Appendix: Proof of Theorem RSC1 and Theorem RSC2}

\begin{proof}
[Proof of Theorem RSC1]

Since the operators $C_\lambda$, $\lambda>0$, defined in \eqref{Clambda_def} are a priori and the operator $C$  is by assumption skew-self-adjoint, we can define the following bounded operators (actually contractions):
\begin{align}\label{eq:Klambda_contraction}
&
K_\lambda:=(I-C_\lambda)^{-1},
&&
\norm{K_\lambda}\le1,
&&
\lambda>0,
\\
&
K:=(I-C)^{-1},
&&
\norm{K}\le1.\notag
\end{align}
Hence,  we can write the resolvent $R_\lambda =(\lambda I-L)^{-1}$ \eqref{resolvent_def} as
\begin{align}
\label{resolvent_master}
R_\lambda
=
(\lambda+S)^{-1/2}
K_\lambda
(\lambda+S)^{-1/2}.
\end{align}

\begin{lemma}
\label{lem:Klambdastcvg}
Assume that the sequence of bounded operators $K_\lambda$ converges to $K$ in the strong operator topology:
\begin{align}
\label{Klambdastcvg}
K_\lambda\tostoptop K,
\qquad \text{as}\qquad\lambda\to0.
\end{align}
Then for any $f\in\Dom(S^{-1/2})=\Ran(S^{1/2})$,  the limits in \eqref{conditionA} hold.
\end{lemma}

\begin{proof}
[Proof of Lemma \ref{lem:Klambdastcvg}]
From the spectral theorem applied to the positive operator $S$, it is obvious that, as $\lambda\to0+$, 
\begin{align}
&
\norm{\lambda^{1/2}(\lambda+S)^{-1/2}}\le1,
&&
\lambda^{1/2}(\lambda+S)^{-1/2}\tostoptop 0,
\notag\\[8pt]
\label{SlambdaS}
&
\norm{S^{1/2}(\lambda+S)^{-1/2}}\le1,
&&
S^{1/2}(\lambda+S)^{-1/2}\tostoptop I.
\end{align}
We can write $f = S^{1/2} g$ with $g\in\cH$.  Now,  using \eqref{resolvent_master}, we get
\begin{align*}
\lambda^{1/2}u_\lambda
&=
\lambda^{1/2}(\lambda+S)^{-1/2}
K_\lambda
(\lambda+S)^{-1/2}S^{1/2} g,
\\[8pt]
S^{1/2}u_\lambda
&=
S^{1/2}(\lambda+S)^{-1/2}
K_\lambda
(\lambda+S)^{-1/2}S^{1/2} g.
\end{align*}
We get
\begin{align*}
S^{1/2}u_\lambda &=
S^{1/2}(\lambda + S)^{1/2} K_\lambda(\lambda+S)^{-1/2}S^{1/2}g
\stackrel{\textrm{$(\ref{SlambdaS})$}}{=}
S^{1/2}(\lambda+S)^{-1/2}K_\lambda(g+o(1))\\
\textrm{By (\ref{Klambdastcvg},\ref{eq:Klambda_contraction})}\qquad
&= S^{1/2}(\lambda+S)^{-1/2}(Kg+o(1))
\stackrel{\textrm{$(\ref{SlambdaS})$}}{=}
Kg+o(1)
\end{align*}
where the notation $o(1)$ is for convergence in norm as $\lambda\to 0$. Verifying the other condition of \eqref{conditionA} is similar.
\end{proof}

In the next lemma, we formulate a sufficient condition for \eqref{Klambdastcvg} to hold. 

\begin{lemma}
\label{lem:strrescvg}
Let $C_n$, $n\in\N$, and $C=C_\infty$ be densely defined closed (possibly unbounded) operators over the Hilbert space $\cH$. Let also $\cC_n$ and $\cC$ be a cores of definition of the operators $C_n$ and $C$, respectively. Assume that some (fixed) $\mu\in\C$ is in the intersection of the resolvent set of all operators $C_n$, $n\le\infty$, and
\begin{align}
\label{unifbound}
\sup_{1\le n\le\infty}\norm{(\mu I - C_n)^{-1}}<\infty,
\end{align}
and for any $h\in\cC$ there exists a sequence $h_n\in\cC_n$ such that the following limits hold
\begin{align}
\label{stcvgoncore}
\lim_{n\to\infty}\norm{h_n-h}=0.
\qquad
\text{and}
\qquad
\lim_{n\to\infty}\norm{C_n h_n -C  h}=0.
\end{align}
Then (i) and (ii) below hold. 
\\
(i)
\begin{align}
\label{strescvg}
(\mu I - C_n)^{-1}
\tostoptop
(\mu I - C)^{-1}.
\end{align}
(ii)
The sequence of operators $C_n$ converges \emph{in the strong graph limit} sense to $C$.
\end{lemma}

\begin{proof}
[Proof of Lemma \ref{lem:strrescvg}]
(i)
Since $\cC$ is a core for the densely defined closed operator $C$ and $\mu$ is in the resolvent set of $C$, the subspace $\wh{\cC}:=\{\wh h=(\mu I  - C) h\,:\,  h\in\cC\}$ is dense in $\cH$. For $\wh h\in\wh\cC$ let $h:= (\mu I - C)^{-1} \wh h\in\cC$ and choose a sequence $h_n\in\cC_n$ for which \eqref{stcvgoncore} holds. Then 
\[
(\mu I-C_n)^{-1} \wh h - (\mu I-C)^{-1} \wh h
=
\left( \mu (\mu I-C_n)^{-1} -I \right) (h-h_n)
+
(\mu I-C_n)^{-1} (C_n h_n -Ch), 
\]
and hence 
\begin{align*}
&
\norm{(\mu I-C_n)^{-1} \wh h - (\mu I-C)^{-1} \wh h}
\\
&\hskip3cm
\le
\left( \abs{\mu}\norm{(\mu I-C_n)^{-1}} +1\right) \norm{h-h_n}
+
\norm{\mu I-C_n)^{-1}} \norm{C_n h_n -Ch}
\to 0. 
\end{align*}
due to \eqref{unifbound} and \eqref{stcvgoncore}. Since this is valid on the \emph{dense} subspace $\wh\cC\subset\cH$, using again \eqref{unifbound}, we conclude \eqref{strescvg}.

\medskip
\noindent
(ii)
The proof of the ``if'' part of Theorem VIII. 26 in \cite{reed_simon_vol1_vol2_75} can be transposed without any essential alteration. 
\end{proof}

\medskip
\noindent
To finish the proof of Theorem RSC1 
first apply Lemma \ref{lem:strrescvg}(i) to 
$C_\lambda$, $\lambda\to0+$, defined in \eqref{Clambda_def}, 
$C$ assumed (essentially) skew self-adjoint, and $\mu=1$.
Note that $\mu=1$ is indeed in the resolvent set of all these operators and, indeed $\sup_{\lambda>0}\norm{(I-C_\lambda)^{-1}}<\infty$ and $\norm{(I-C)^{-1}}<\infty$, as required in \eqref{unifbound}, since the operators $C_\lambda$ are bounded and skew-self-adjoint and the operator $C$ is assumed to be essentially skew-self-adjoint. From Lemma \ref{lem:strrescvg}(i) it follows that that \eqref{Klambdastcvg} holds. Finally, quoting Lemma \ref{lem:Klambdastcvg} we conclude the proof of Theorem RSC1. 
\end{proof}

\begin{proof}[Proof of Theorem RSC2]
From $0\le T\le cD$ \eqref{Ddominates} it follows that
\begin{align}
\label{Ddominates2}
0\le D \le S \le (1+c) D
\end{align}
Let
\begin{align}
\notag
&
V_\lambda:=
(\lambda I+D)^{1/2}
(\lambda I+S)^{-1/2},
&&
V=V_0:=
D^{1/2}
S^{-1/2}.
\end{align}
The operator $V$ is a priori defined on $\Dom(S^{-1/2})=\Ran(S^{1/2})$, but as we see next, it extends by continuity to a bounded and invertible linear operator defined on the whole space $\cH$. Due to \eqref{Ddominates2} the following bounds hold uniformly for $\lambda\ge0$:
\begin{align}
\notag
&
\norm{V_\lambda}
=
\norm{V_\lambda^*}
\le 1,
&&
\norm{V_\lambda^{-1}}
=
\norm{(V_\lambda^{-1})^*}
\le
\sqrt{1+c}.
\end{align}
Let us show that bound on $\norm{V_\lambda}$, the bound on
$\norm{V_\lambda^{-1}}$ is similar. We write
\begin{align*}
\norm{V_\lambda\varphi}^2
&=\langle(\lambda I+D)^{1/2}(\lambda I+S)^{-1/2}\varphi,
(\lambda I+D)^{1/2}(\lambda I+S)^{-1/2}\varphi\rangle \\
& =
\langle(\lambda I+S)^{-1/2}\varphi,
(\lambda I+D)(\lambda I+S)^{-1/2}\varphi\rangle\\
&\le 
\langle(\lambda I+S)^{-1/2}\varphi,
(\lambda I+S)(\lambda I+S)^{-1/2}\varphi\rangle
=\norm{\varphi}^2.
\end{align*}
From here, first of all,  it follows that
\begin{align}
\notag
\Dom(S^{-1/2})=\Dom(D^{-1/2}),
\end{align}
and thus the $\cH_{-1}$-conditions $f\in\Dom(S^{-1/2})$, respectively, $f\in\Dom(D^{-1/2})$ in Theorem RSC1, respectively, Theorem RSC2, are actually the same. It is also easy to see that for any $\varphi\in\cH$
\[
\lim_{\lambda\to0}
V_\lambda\varphi= V\varphi
\qquad
\text{and}
\qquad
\lim_{\lambda\to0}
V_\lambda^{-1}\varphi= V^{-1}\varphi. 
\]
That is, $V_\lambda\tostoptop V$ and $V^{-1}_\lambda\tostoptop V^{-1}$, as $\lambda\to0$, where  $\tostoptop$ stands for convergence in the strong operator topology.

Next write the operators $C_\lambda$ and $C$ from Theorem RSC1, as
\begin{align}
\notag
&
C_\lambda
=
V^*_{\lambda} B_\lambda V_{\lambda},
&&
C
=
V^* B V.
\end{align}
Now, from the fact that $V_\lambda$ and $V_\lambda^{-1}$ are all bounded, uniformly in $\lambda\ge0$, it readily follows that:
(a) one can use $\cC=V^{-1}\cB$ as a core for the operator $C$;
(b) $C$ is essentially skew-self-adjoint on $\cC$ if so was $B$ on $\cB$; and
(c) the limit \eqref{Clambdalimit} follows from \eqref{Blambdalimit} by straightforward manipulations. Indeed, for $\psi\in\cC$ define $\varphi:=V\psi\in\cB$ and let $\varphi_\lambda \in\cH$ be such that the limits in \eqref{Blambdalimit} hold. Define $\psi_\lambda:=V_\lambda^{-1}\varphi_\lambda$. Then the limits in \eqref{Clambdalimit} clearly hold:
\begin{align*}
\norm{\psi_\lambda-\psi}
&=
\norm{V_\lambda^{-1}\varphi_\lambda - V^{-1}\varphi}
\le
\norm{V_\lambda^{-1}}\norm{\varphi_\lambda-\varphi}
+
\norm{V^{-1}_\lambda\varphi-V^{-1}\varphi}
\to0, 
\\[10pt]
\norm{C_\lambda\psi_\lambda-C\psi}
&=
\norm{V^{*}_\lambda B_\lambda\varphi_\lambda - V^{*}B\varphi}
\le
\norm{V_\lambda^{*}}\norm{B_\lambda\varphi_\lambda-B\varphi}
+
\norm{V^{*}_\lambda B\varphi-V^{*} B\varphi}
\to0.\qedhere
\end{align*}
\end{proof}

\bigskip\bigskip

\noindent
{\bf Acknowledgements:}
We thank an anonymous reviewer for thorough criticism, and in particular for the recommendation to formulate a more generally valid version of the main result. 
The research of BT is partially supported by OTKA (HU) grant K 100473 and by EPSRC (UK) Fellowship, grant no. EP/P003656/1. 
The research of GK is partially supported by the Israel Science Foundation and the Jesselson Foundation. 
Both authors acknowledge mobility support by The Leverhulme Trust (UK) through the International Network ``Laplacians, Random Walks, Quantum Spin Systems''.

\vskip2cm

\hbox{
\vbox{\hsize=7cm\noindent
{\sc Gady Kozma}
\\
Department of Mathematics
\\
The Weizmann Institute of Science
\\
POB 26, Rehovot, 76100
\\
Israel
\\
email: {\tt gady.kozma@weizmann.ac.il}
}
\hskip2cm
\vbox{\hsize=7cm\noindent
{\sc B\'alint T\'oth}
\\
School of Mathematics
\\
University of Bristol
\\
Bristol, BS8 1TW
\\
United Kingdom
\\
email: {\tt balint.toth@bristol.ac.uk}
}
}

\end{document}